\documentclass[aap,seceqn,MSNbibl,citesort,dvips]{arximspdf}
\usepackage{graphicx}

\doi{10.1214/09-AAP667}
\volume{20}
\issue{5}
\pubyear{2010}
\firstpage{1607}
\lastpage{1637}

\makeatletter

\newproclaim{assumption}{Assumption}[section]
\newproclaim{remark}[assumption]{Remark}
\newtheorem{lemma}[assumption]{Lemma}
\newproclaim{definition}[assumption]{Definition}
\newtheorem{theorem}[assumption]{Theorem}
\newtheorem{proposition}[assumption]{Proposition}
\newtheorem{corollary}[assumption]{Corollary}

\makeatother

\begin{document}
\begin{frontmatter}

\title{Numerical method for optimal stopping of piecewise deterministic
Markov processes\thanksref{T1}}
\runtitle{Optimal stopping for PDMP}

\thankstext{T1}{Supported by ARPEGE program of the French National
Agency of Research (ANR),
project ``FAUTOCOES,'' number ANR-09-SEGI-004.}

\begin{aug}
\author[A]{\fnms{Beno\^{\i}te} \snm{de~Saporta}\ead[label=e2]{saporta@math.u-bordeaux1.fr}},
\author[A]{\fnms{Fran\c{c}ois} \snm{Dufour}\corref{}\ead[label=e1]{dufour@math.u-bordeaux1.fr}} and
\author[A]{\fnms{Karen} \snm{Gonzalez}\ead[label=e3]{gonzalez@math.u-bordeaux1.fr}}
\runauthor{B. de Saporta, F. Dufour and K. Gonzalez}
\affiliation{Universit\'{e} de Bordeaux}
\address[A]{B. de Saporta\\
F. Dufour\\
K. Gonzalez\\
Universit\'{e} de Bordeaux\\
IMB, 351 cours de la Lib\'{e}ration\\
F33405 Talence Cedex\\
France\\
\printead{e1}\\
\phantom{E-mail: }\printead*{e2}\\
\phantom{E-mail: }\printead*{e3}} 
\end{aug}

\pdfauthor{Benoite de Saporta, Francois Dufour, Karen Gonzalez}

\received{\smonth{3} \syear{2009}}
\revised{\smonth{9} \syear{2009}}

%
\begin{abstract}
We propose a numerical method to approximate the value function for the
optimal stopping problem of a piecewise deterministic Markov process
(PDMP). Our approach is based on quantization of the post jump
location---inter-arrival time Markov chain naturally embedded in the
PDMP, and
path-adapted time discretization grids. It allows us to derive bounds
for the convergence rate of the algorithm and to provide a computable
$\epsilon$-optimal stopping time. The paper is illustrated by a
numerical example.
\end{abstract}

%

\end{frontmatter}

\section{Introduction}
The aim of this paper is to propose a computational method for optimal
stopping of a piecewise deterministic Markov process $\{X(t)\}$ by
using a quantization technique for an underlying discrete-time Markov
chain related to the continuous-time process $\{X(t)\}$ and
path-adapted time discretization grids.

Piecewise-deterministic Markov processes (PDMPs) have been introduced
in the literature by Davis \cite{davis93} as a general class of
stochastic models.
PDMPs are a family of Markov processes involving deterministic motion
punctuated by random jumps.
The motion of the PDMP $\{X(t)\}$ depends on three local
characteristics, namely the flow $\phi$, the jump rate $\lambda$ and the
transition measure $Q$, which specifies the post-jump location.
Starting from $x$ the motion of the process follows the flow $\phi
(x,t)$ until the first jump time $T_1$ which occurs either
spontaneously in a Poisson-like fashion with rate $\lambda(\phi(x,t))$
or when the flow $\phi(x,t)$ hits the boundary of the state-space.
In either case the location of the process at the jump time $T_1\dvtx
X(T_{1})=Z_{1}$ is selected by the transition measure $Q(\phi
(x,T_1),\cdot)$.
Starting from $Z_{1}$, we now select the next interjump time
$T_{2}-T_{1}$ and postjump location $X(T_{2})=Z_{2}$.
This gives a piecewise deterministic trajectory for $\{X(t)\}$ with
jump times $\{T_{k}\}$ and post jump locations $\{Z_{k}\}$ which
follows the flow $\phi$ between two jumps.
A suitable choice of the state space and the local characteristics
$\phi
$, $\lambda$ and $Q$ provide stochastic
models covering a great number of problems of operations research \cite
{davis93}.

Optimal stopping problems have been studied for PDMPs in \cite
{costa88,costa00,davis93,gatarek91,gugerli86,lenhart85}. In \cite
{gugerli86} the author defines an operator related to the first jump
time of\vadjust{\goodbreak} the process and shows that the value function of the optimal
stopping problem is a fixed
point for this operator.
The basic assumption in this case is that the final cost function is
continuous along trajectories, and it is shown that the
value function will also have this property.
In \cite{gatarek91,lenhart85} the authors adopt some stronger
continuity assumptions and boundary conditions to show that the value
function of the optimal stopping
problem satisfies some variational inequalities related to
integro-differential equations.
In \cite{davis93}, Davis assumes that the value function is bounded and
locally Lipschitz along trajectories to show that the variational
inequalities are necessary and sufficient to
characterize the value function of the optimal stopping problem.
In~\cite{costa00}, the authors weakened the continuity assumptions of
\cite{davis93,gatarek91,lenhart85}.
A~paper related to our work is \cite{costa88} by Costa and Davis. It is
the only one presenting a
computational technique for solving the optimal stopping problem for a
PDMP based on a discretization of the state space similar to the one
proposed by Kushner in \cite{kushner77}.
In particular, the authors in \cite{costa88} derive a convergence
result for the approximation scheme but no estimation of the rate of
convergence is derived.

Quantization methods have been developed recently in numerical
probability, nonlinear filtering or optimal stochastic control with
applications in finance
\cite{bally03,bally05,pages98,pages05,pages04b,pages04}.
More specifically, powerful and interesting methods have been developed
in \mbox{\cite{bally03,bally05,pages04}} for computing the Snell-envelope
associated to discrete-time Markov chains and diffusion processes.
Roughly speaking, the approach developed in \cite
{bally03,bally05,pages04} for studying the optimal stopping problem for
a continuous-time diffusion process $\{Y(t)\}$ is based on a
time-discretization scheme to obtain a discrete-time Markov chain $\{
\overline{Y}_{k}\}$. It is shown that the original continuous-time
optimization problem can be converted to an auxiliary optimal stopping
problem associated with the discrete-time Markov chain $\{\overline
{Y}_{k}\}$. Under some suitable assumptions, a~rate of convergence of
the auxiliary value function to the original one can be derived.
Then, in order to address the optimal stopping problem of the
discrete-time Markov chain, a twofold computational method is proposed.
The first step consists in approximating the Markov chain by a
quantized process. There exists an extensive literature on quantization
methods for random variables and processes. We do not pretend to
present here an exhaustive panorama of these methods. However, the
interested reader may, for instance, consult
\cite{gray98,pages98,pages04} and the references therein.
The second step is to approximate the conditional expectations which
are used to compute the backward dynamic programming formula by the
conditional expectation related to the quantized process. This
procedure leads to a tractable formula called a \textit{quantization
tree algorithm}
(see Proposition 4 in \cite{bally03} or Section~4.1 in \cite{pages04}).
Providing the cost function and the Markov kernel are Lipschitz, some
bounds and rates of convergence are obtained
(see, e.g., Section~2.2.2 in \cite{bally03}).

As regards PDMPs, it was shown in \cite{gugerli86} that the value
function of the optimal stopping problem can be calculated by iterating
a functional operator, labeled $L$ [see (\ref{expL}) for its
definition], which is related to a continuous-time maximization and a
discrete-time dynamic programming formula. Thus, in order to
approximate the value function of the optimal stopping problem of a
PDMP $\{X(t)\}$, a natural approach would have been to follow the same
lines as in
\cite{bally03,bally05,pages04}.
However, their method cannot be directly applied to our problem for two
main reasons related to the specificities of PDMPs.

First, PDMPs are in essence discontinuous at random times. Therefore,
as pointed out in \cite{gugerli86},
it will be problematic to convert the original optimization problem
into an optimal stopping problem associated to a time discretization
of $\{X(t)\}$ with nice convergence properties.
In particular, it appears ill-advised to propose as in \cite{bally03} a
fixed-step time-discretization scheme $\{X(k\Delta)\}$ of the original
process $\{X(t)\}$.
Besides, another important intricacy concerns the transition semigroup
$\{P_{t}\}_{t\in\mathbb{R}_{+}}$ of $\{X(t)\}$. On the one hand, it cannot
be explicitly calculated from the local characteristics $(\phi,\lambda
,Q)$ of the PDMP (see \cite{costa08,dufour99}).
Consequently, it will be complicated to express the Markov kernel
$P_{\Delta}$ associated with the Markov chain $\{X(k\Delta)\}$. On the
other hand, the Markov chain $\{X(k\Delta)\}$ is, in general, not even
a Feller chain (see
\cite{davis93}, pages 76 and 77), and therefore it will be hard to
ensure it is $K$-Lipschitz (see Definition 1 in \cite{bally03}).

Second, the other main difference stems from the fact that the function
appearing in the backward dynamic programming formula
associated with $L$ and the reward function $g$ is not continuous even
if some strong regularity assumptions are made on $g$.
Consequently, the approach developed in \cite{bally03,bally05,pages04}
has to be refined since it can only handle conditional expectations of
Lipschitz-continuous functions.

However, by using the special structure of PDMPs, we are able to
overcome both these obstacles. Indeed, associated to the PDMP $\{X(t)\}
$, there exists a natural embedded discrete-time Markov chain $\{\Theta
_{k}\}$ with $\Theta_{k}=(Z_{k},S_{k})$ where
$S_{k}$ is given by the inter-arrival time $T_{k}-T_{k-1}$. The main
operator $L$ can be expressed using the chain $\{\Theta_{k}\}$ and a
continuous-time maximization. We first convert the continuous-time
maximization of operator $L$ into a discrete-time maximization by using
a path-dependent time-discretization scheme. This enables us to
approximate the value function by the solution of a backward dynamic
programming equation in discrete-time involving conditional expectation
of the Markov chain~$\{\Theta_{k}\}$.
Then, a~natural approximation of this optimization problem is obtained
by replacing $\{\Theta_{k}\}$ by its quantized approximation.
It must be pointed out that this optimization problem is related to the
calculation of conditional expectations of indicator functions of the
Markov chain $\{\Theta_{k}\}$.
As said above, it is not straightforward to obtain convergence results
as in \cite{bally03,bally05,pages04}.
We deal successfully with indicator functions by showing that the event
on which the discontinuity actually occurs is of small enough
probability. This enables us to provide a rate of convergence for the
approximation scheme.

In addition, and more importantly, this numerical approximation scheme
enables us to propose a computable stopping rule which also is an
$\epsilon$-optimal stopping time of the original stopping problem.
Indeed, for any $\epsilon>0$ one can construct a stopping time,
labeled $\tau$, such that
\[
V(x)-\epsilon\leq\mathbf{E}_{x} [g(X(\tau)) ]\leq V(x),
\]
where $V(x)$ is the optimal value function associated to the original
stopping problem.
Our computational approach is attractive in the sense that it does not
require any additional calculations.
Moreover, we can characterize how far it is from optimal in terms of
the value function.
In \cite{bally03}, Section 2.2.3, Proposition 6, another criteria for
the approximation of the optimal stopping time has been proposed.
In the context of PDMPs, it must be noticed that an optimal stopping
time does not generally exist as shown in \cite{gugerli86}, Section~2.

An additional result extends Theorem 1 of Gugerli \cite{gugerli86} by
showing that the iteration of operator $L$ provides
a sequence of random variables which corresponds to a \textit
{quasi}-Snell envelope associated
with the reward process $ \{g(X(t)) \}_{t\in\mathbb{R}_{+}}$ where the
horizon time is random and given by the jump times $(T_{n})_{n\in\{
0,\ldots,N\}}$ of the process $ \{X(t) \}_{t\in\mathbb{R}_{+}}$.

The paper is organized as follows. In Section \ref{section_defi} we
give a precise definition of PDMPs and state our notation and
assumptions. In Section \ref{sectionoptstop}, we state the optimal
stopping problem, recall and refine some results from
\cite{gugerli86}. In Section \ref{section_Vchap}, we build an
approximation of the value function. In Section \ref{section_erreur},
we compute the error between the approximate value function and the
real one. In Section \ref{sectiontpsarret} we propose a computable
$\epsilon$-optimal stopping time and evaluate its sharpness. Finally in
Section \ref{section_appli} we present a numerical example. Technical
results are postponed to the \hyperref[app]{Appendix}.

\section{Definitions and assumptions}
\label{section_defi}
We first give a precise definition of a piecewise deterministic Markov process.
Some general assumptions are presented in the second part of this section.
Let us introduce first some standard notation. Let $M$ be a metric
space. $\mathbf{B}(M)$ is the set of real-valued, bounded, measurable
functions defined on $M$. The Borel $\sigma$-field of $M$ is denoted by
$\mathcal{B}(M)$. Let $Q$ be a Markov kernel on $(M,\mathcal{B}(M))$
and $w\in\mathbf{B}(M)$, $Qw(x)=\int_{M} w(y) Q(x,dy)$ for $x\in M$.
For $(a,b)\in\mathbb{R}^2$, $a\wedge b= \min(a,b)$ and $a\vee b=
\max(a,b)$.

\subsection{Definition of a PDMP}
Let $E$ be an open subset of $\mathbb{R}^n$, $\partial E$ its boundary and
$\overline{E}$ its closure.
A PDMP is determined by its local characteristics $(\phi,\lambda,Q)$ where:

$\bullet$ The flow $\phi\dvtx \mathbb{R}^{n} \times\mathbb{R}\to
\mathbb{R}^{n}$ is a
one-parameter group of homeomorphisms: $\phi$ is continuous,
$\phi(\cdot, t)$ is an homeomorphism for each $t\in\mathbb{R}$ satisfying
$\phi(\cdot, t+s)=\phi(\phi(\cdot, s),t))$.

For all $x$ in $E$, let us denote
\[
t^{*}(x)\doteq\inf\{t>0\dvtx\phi(x,t)\in\partial E \}
\]
with the convention $\inf\varnothing= \infty$.

$\bullet$ The jump rate $\lambda\dvtx \overline{E} \to\mathbb{R}_{+}$
is assumed
to be a measurable function satisfying
\[
(\forall x \in E),\qquad (\exists\varepsilon>0) \qquad\mbox{such that }
\int_0^\varepsilon\lambda(\phi(x,s))\,ds< \infty.
\]

$\bullet$ $Q$ is a Markov kernel on $(\overline{E},\mathcal
{B}(\overline{E}))$ satisfying the following property:
\[
(\forall x\in\overline{E}),\qquad Q(x,E-\{x\})=1.
\]
From these characteristics, it can be shown \cite{davis93}, pages
62--66, that there exists a filtered probability space $(\Omega
,\mathcal
{F},\{ \mathcal{F}_{t} \}, \{ \mathbf{P}_{x} \}_{x\in E})$
such that the motion of the process $\{X(t)\}$ starting from a point
$x\in E$ may be constructed as follows.
Take a random variable $T_1$ such that
\[
\mathbf{P}_x(T_1>t) \doteq
\cases{
e^{-\Lambda(x,t)}, &\quad for $t<t^{*}(x)$,\cr
0, &\quad for $t\geq t^{*}(x)$,}
\]
where for $x\in E$ and $t\in[0,t^{*}(x)]$
\[
\Lambda(x,t) \doteq\int_0^t\lambda(\phi(x,s))\,ds.
\]
If $T_1$ generated according to the above probability is equal to
infinity, then for $t\in\mathbb{R}_+$, $X(t)=\phi(x,t)$.
Otherwise select independently an $E$-valued random variable (labelled
$Z_{1}$) having distribution $Q(\phi(x,T_1),\cdot)$,
namely $\mathbf{P}_x(Z_{1}\in A)=Q(\phi(x,T_1),A)$ for any $A\in
\mathcal{B}(\overline{E})$.
The trajectory of $\{X(t)\}$ starting at $x$, for $t\leq T_1$, is
given by
\[
X(t) \doteq
\cases{
\phi(x,t), &\quad for $t<T_1$, \cr
Z_1, &\quad for $t=T_1$.}
\]
Starting from $X(T_1)=Z_1$, we now select the next inter-jump time $T_2-T_1$
and post-jump location $X(T_2)=Z_2$ is a similar way.

This gives a strong Markov process $\{X(t)\}$ with jump times $\{
T_k\}_{k \in\mathbb{N}}$ (where $T_{0}= 0$).
Associated with $\{X(t)\}$, there exists a discrete time process $ (
\Theta_{n} )_{n\in\mathbb{N}}$ defined by
$\Theta_{n}=(Z_{n},S_{n})$ with $Z_{n}=X(T_{n})$ and
$S_{n}=T_{n}-T_{n-1}$ for $n\geq1$ and $S_{0}=0$.
Clearly, the process $(\Theta_{n})_{n\in\mathbb{N}}$ is a Markov chain.

We introduce a standard assumption (see, e.g., equations (24.4) or
(24.8) in \cite{davis93}).
\begin{assumption}
\label{A1}
For all $(x,t)\in E\times\mathbb{R}_{+}$, $\mathbf{E}_{x} [ \sum
_{k} \mathbf{1}
_{\{T_{k}\leq t\}} ] < \infty$.
\end{assumption}

In particular, it implies that $T_k \to\infty$ as $k \to\infty$.

For $n\in\mathbb{N}$, let $\mathcal{M}_{n}$ be the family of all $\{
\mathcal
{F}_{t} \}$-stopping times which are dominated by $T_{n}$, and
for $n<p$, let $\mathcal{M}_{n,p}$ be the family of all $\{ \mathcal
{F}_{t} \}$-stopping times $\nu$ satisfying
$T_{n}\leq\nu\leq T_{p}$. Let $\mathbf{B}^{c}$ denote the set of all
real-valued, bounded, measurable functions, $w$ defined on $\overline
{E}$ and continuous along trajectories up to the jump time horizon: for
any $x\in E$, $w(\phi(x,\cdot))$ is continuous on $ [0,t^{*}(x)]$. Let
$\mathbf{L}^{c}$ be the set of all real-valued, bounded, measurable
functions, $w$ defined on $\overline{E}$ and
Lipschitz along trajectories:
\begin{enumerate}
\item there exists $[ w ]_{1}\in\mathbb{R}_{+}$ such that for any
$(x,y)\in
E^{2}$, $u\in[0,t^{*}(x)\wedge t^{*}(y)]$, one has
\[
| w(\phi(x,u))-w(\phi(y,u)) | \leq[ w ]_{1} |x-y|;
\]
\item there exists $[ w ]_{2}\in\mathbb{R}_{+}$ such that for any
$x\in E$,
and $(t,s)\in[0,t^{*}(x)]^{2}$, one has
\[
| w(\phi(x,t))-w(\phi(x,s)) | \leq[ w ]_{2} |t-s|;
\]
\item there exists $[ w ]_{*}\in\mathbb{R}_{+}$ such that for any
$(x,y)\in
E^{2}$, one has
\[
| w(\phi(x,t^{*}(x)))-w(\phi(y,t^{*}(y))) | \leq[ w ]_{*} |x-y|.
\]
\end{enumerate}

In the sequel, for any function $f$ in $\mathbf{B}^{c}$, we denote by
$C_f$ its bound
\[
C_f={\sup_{x\in E}}|f(x)|,
\]
and for any Lipschitz-continuous function $f$ in $\mathbf{B}(E)$ or
$\mathbf{B}(\overline{E})$, we denote by $[f]$ its Lipschitz constant
\[
[ f ]=\sup_{x\neq y\in E}\frac{ |f(x)-f(y) |}{|x-y|}.
\]
\begin{remark}\label{remarque_lip}
$\mathbf{L}^{c}$ is a subset of $\mathbf{B}^{c}$ and any function in
$\mathbf{L}^{c}$ is Lipschitz on $\overline{E}$ with $[ w ]\leq[ w ]_1$.
\end{remark}

Finally, as a convenient abbreviation, we set for any $x\in\overline
{E}$, $\lambda Qw(x) = \lambda(x) Qw(x)$.

\subsection{Assumptions}
The following assumptions will be in force throughtout.

\begin{assumption}
\label{H2}
The jump rate $\lambda$ is bounded and there exists $[ \lambda
]_{1}\in\mathbb{R}_{+}$ such that for any
$(x,y)\in E^{2}$, $u\in[0,t^{*}(x)\wedge t^{*}(y)[$,
\[
| \lambda(\phi(x,u))-\lambda(\phi(y,u)) | \leq[ \lambda ]_{1} |x-y|.
\]
\end{assumption}
\begin{assumption}
\label{H3}
The exit time $t^{*}$ is bounded and Lipschitz-continuous on~$E$.
\end{assumption}
\begin{assumption}
\label{H5}
The Markov kernel $Q$ is Lipschitz in the following sense: there exists
$[ Q ]\in\mathbb{R}_{+}$ such that for any function $w\in\mathbf{L}^{c}$
the following two conditions are satisfied:
\begin{enumerate}
\item for any $(x,y)\in E^{2}$, $u\in[0,t^{*}(x)\wedge t^{*}(y)]$,
one has
\[
| Qw(\phi(x,u))-Qw(\phi(y,u)) | \leq[ Q ] [ w ]_{1} |x-y|;
\]
\item for any $(x,y)\in E^{2}$, one has
\[
| Qw(\phi(x,t^{*}(x)))-Qw(\phi(y,t^{*}(y))) | \leq[ Q ] [ w ]_{*} |x-y|.
\]
\end{enumerate}
\end{assumption}

The reward function $g$ associated with the optimal stopping problem
satisfies the following hypothesis.
\begin{assumption}
\label{H1}
$g$ is in $\mathbf{L}^{c}$.
\end{assumption}

\section{Optimal stopping problem}
\label{sectionoptstop}
From now on, assume that the distribution of $X(0)$ is given by $\delta
_{x_{0}}$ for a fixed state $x_{0}\in E$.
Let us consider the following optimal stopping problem for a fixed
integer $N$:
%
%
\begin{equation}\label{DefOpt}
\sup_{\tau\in\mathcal{M}_{N}}\mathbb{E}_{x_{0}}[g(X(\tau))].
\end{equation}
This problem has been studied by Gugerli   \cite{gugerli86}.

Note that Assumption \ref{H2} yields $\Lambda(x,t)<\infty$ for all
$x,t$. Hence, for all $x$ in~$E$,
the jump time horizon $s^{*}(x)$ defined in \cite{gugerli86} by $t^*(x)
\wedge\inf\{t\geq0, {e}^{-\Lambda(x,t)}=0\}$ is equal to the exit
time $t^{*}(x)$.
Therefore, operators $H\dvtx\mathbf{B}(\overline{E}) \to\mathbf
{B}(E\times
\mathbb{R}_{+})$, $I\dvtx\mathbf{B}(E) \to\mathbf{B}(E\times\mathbb{R}_{+})$,
$J\dvtx\mathbf{B}(E) \times\mathbf{B}(\overline{E}) \to\mathbf
{B}(E\times
\mathbb{R}_{+})$, $K\dvtx\mathbf{B}(E) \to\mathbf{B}(E)$
and $L\dvtx\mathbf{B}(E) \times\mathbf{B}^{c} \to\mathbf{B}^{c}$
introduced by Gugerli   (\cite{gugerli86}, Section 2) reduce to
%
%
\begin{eqnarray}
Hf(x,t)&=&f \bigl(\phi\bigl(x,t\wedge t^{*}(x)\bigr) \bigr){e}^{-\Lambda(x,t\wedge
t^{*}(x))},\nonumber\\
Iw(x,t)&=&\int_{0}^{t\wedge t^{*}(x)} \lambda Qw(\phi(x,s))
{e}^{-\Lambda(x,s)}\,ds,\nonumber\\
\label{defJ}
J(w,f)(x,t)&=&Iw(x,t)+Hf(x,t),
\\
\label{defK}
Kw(x)&=&\int_{0}^{t^{*}(x)} \lambda Qw(\phi(x,s))
{e}^{-\Lambda(x,s)}\,ds\nonumber\\[-8pt]\\[-8pt]
&&{}+Qw(\phi(x,t^{*}(x))) {e}^{-\Lambda(x,t^{*}(x))},\nonumber
\\
L(w,h)(x) & = & \sup_{t\geq0}J(w,h)(x,t)\vee Kw(x). \nonumber
\end{eqnarray}

It is easy to derive a probabilistic interpretation of operators $H$,
$I$, $K$ and $L$ in terms of the embedded Markov chain
$ ( Z_{n},S_{n} )_{n\in\mathbb{N}}$.
\begin{lemma}
\label{Iproba}
For all $x\in E$, $w\in\mathbf{B}(E)$, $f\in\mathbf{B}(\overline{E})$,
$h\in\mathbf{B}^{c}$ and $t\geq0$, one has
%
%
\begin{eqnarray}
Hf(x,t) & = & f \bigl(\phi\bigl(x,t\wedge t^{*}(x)\bigr) \bigr) \mathbf{P}_{x}
\bigl(
S_{1}\geq t\wedge t^{*}(x) \bigr), \nonumber\\
Iw(x,t)& = & \mathbf{E}_{x} \bigl[w(Z_{1}) \mathbf{1}_{ \{S_{1}< t
\wedge
t^{*}(x) \} } \bigr], \nonumber\\
\label{Kproba}
Kw(x)& = & \mathbf{E}_{x} [w(Z_{1}) ], \\
\label{expL}
L(w,h)(x) &=&\sup_{u\leq t^{*}(x)} \bigl\{ \mathbf{E}_{x} \bigl[w(Z_{1})
\mathbf{1}_{ \{ S_{1}<u \} } \bigr]
+ h (\phi(x,u) ) \mathbf{P}_{x} ( S_{1}\geq u )
\bigr\}\nonumber\\[-8pt]\\[-8pt]
&&{}
\vee\mathbf{E}_{x} [w(Z_1) ].\nonumber
\end{eqnarray}
\end{lemma}

For a reward function $g \in\mathbf{B}^{c}$, it has been shown in
\cite
{gugerli86} that the value function can be recursively constructed by
the following procedure:
\[
\sup_{\tau\in\mathcal{M}_{N}}\mathbb{E}_{x_{0}}[g(X(\tau))] =
v_{0}(x_{0})
\]
with
\[
\cases{
v_{N}=g, \cr
v_{k}=L(v_{k+1},g), &\quad for $k\leq N-1$.}
\]
\begin{definition}
Introduce the random variables $(V_{n})_{n\in\{0,\ldots,N\}}$ by
\[
V_{n}=v_{n}(Z_{n})
\]
or equivalently
%
%
\begin{eqnarray}\label{snell}
V_{n} & = & \sup_{u\leq t^{*}(Z_{n})} \bigl\{ \mathbf{E} \bigl[
v_{n+1}(Z_{n+1}) \mathbf{1}_{\{ S_{n+1} <u \} }
+g (\phi(Z_{n},u) )\mathbf{1}_{\{ S_{n+1}\geq u\} } | Z_{n}
\bigr] \bigr\} \nonumber\\[-8pt]\\[-8pt]
& &{} \vee\mathbf{E} [v_{n+1}(Z_{n+1}) | Z_{n} ].\nonumber
\end{eqnarray}
\end{definition}

The following result shows that the sequence $(V_{n})_{n\in\{0,\ldots
,N\}}$ corresponds to a \textit{quasi}-Snell envelope associated
with the reward process $ \{g(X(t)) \}_{t\in\mathbb{R}_{+}}$ where the
horizon time is random and given by the jump times $(T_{n})_{n\in\{
0,\ldots,N\}}$ of the process $ \{X(t) \}_{t\in\mathbb{R}_{+}}$:
\begin{theorem}
Consider an integer $n< N$. Then
\[
V_{n} = \sup_{\nu\in\mathcal{M}_{n,N}} \mathbf
{E}_{x_0}[g(X(\nu
))| \mathcal{F}_{T_{n}}].
\]
\end{theorem}
\begin{pf}
Let $\nu\in\mathcal{M}_{n,N}$.
According to Proposition \ref{Aprop1} and Corollary \ref{Acoro1} in
Appendix \ref{app_B}, there exists
$\widehat{\nu} \dvtx E\times(\mathbb{R}_{+}\times E)^{n} \times
\Omega
\to
\mathbb{R}_{+}$ such that for all
$(z_{0},\gamma)\in E\times(\mathbb{R}_{+}\times E)^{n}$ the mapping
$\widehat
{\nu}(z_{0},\gamma) \dvtx\Omega\to\mathbb{R}_{+}$
is an $\{\mathcal{F}_{t}\}_{t\in\mathbb{R}_{+}}$-stopping time
satisfying $\widehat{\nu}(z_{0},\gamma) \leq T_{N-n}$, and $\nu
=T_{n} +
\widehat{\nu}(Z_{0},\Gamma_{n},\theta_{T_{n}})$, where $\Gamma
_{n}=
(S_{1},Z_{1},\ldots,S_{n},Z_{n} )$
and $\theta$ is the shift operator. For $(z_{0},\gamma)\in E\times
(\mathbb{R}
_{+}\times E)^{n}$ define $\mathcal{W} \dvtx E\times(\mathbb
{R}_{+}\times
E)^{n} \to\mathbb{R}$ by
\[
\mathcal{W}(z_{0},\gamma) = \mathbf{E}_{z_{n}}[g(X(\widehat{\nu
}(z_{0},\gamma)))]\leq
\sup_{\tau\in\mathcal{M}_{N-n}} \mathbf{E}_{Z_{n}}[g(X(\tau))],
\]
where $\gamma=(s_{1},z_{1},\ldots,s_{n},z_{n})$.
Hence, the strong Markov property of the process $\{X(t)\}$ yields
\[
\mathbf{E}_{x_0}[g(X(\nu))|\mathcal{F}_{T_{n}}]
= \mathbf{E}_{x_0}\bigl[g\bigl(X\bigl(T_{n}+\widehat{\nu}(Z_{0},\Gamma_{n},\theta
_{T_{n}})\bigr)\bigr)|\mathcal{F}_{T_{n}}\bigr]
= \mathcal{W}(Z_{0},\Gamma_{n}).
\]
Consequently, one has
\[
\mathbf{E}_{x_0}[g(X(\nu))|\mathcal{F}_{T_{n}}] \leq\sup_{\tau
\in
\mathcal{M}_{N-n}}
\mathbf{E}_{Z_{n}}[g(X(\tau))]
\]
and, therefore, one has
%
%
\begin{equation}\label{equal1}
\sup_{\nu\in\mathcal{M}_{n,N}} \mathbf{E}_{x_0}[g(X(\nu
))|\mathcal
{F}_{T_{n}}] \leq
\sup_{\tau\in\mathcal{M}_{N-n}} \mathbf{E}_{Z_{n}}[g(X(\tau))].
\end{equation}

Conversely, consider $\tau\in\mathcal{M}_{N-n}$.
It is easy to show that $T_{n}+\tau\circ\theta_{T_{n}} \in\mathcal
{M}_{n,N}$.
The strong Markov property of the process $\{X(t)\}$ again yields
\[
\mathbf{E}_{Z_{n}}[g(X(\tau))]
= \mathbf{E}_{x_0}\bigl[g\bigl(X(T_{n}+\tau\circ\theta_{T_{n}} )\bigr)|\mathcal
{F}_{T_{n}}\bigr]
\leq\sup_{\nu\in\mathcal{M}_{n,N}} \mathbf{E}_{x_0}[g(X(\nu
))|\mathcal{F}_{T_{n}}]
\]
and hence we obtain
%
%
\begin{equation}\label{equal2}
\sup_{\tau\in\mathcal{M}_{N-n}} \mathbf{E}_{Z_{n}}[g(X(\tau))]
\leq
\sup_{\nu\in\mathcal{M}_{n,N}} \mathbf{E}_{x_0}[g(X(\nu
))|\mathcal
{F}_{T_{n}}].
\end{equation}
Combining equations (\ref{equal1}) and (\ref{equal2}), one has
\[
\sup_{\tau\in\mathcal{M}_{N-n}} \mathbf{E}_{Z_{n}}[g(X(\tau))]
=
\sup_{\nu\in\mathcal{M}_{n,N}} \mathbf{E}_{x_0}[g(X(\nu
))|\mathcal
{F}_{T_{n}}].
\]
Finally, it is proved in \cite{gugerli86}, Theorem 1, that
$v_{n}(x)=\sup_{\tau\in\mathcal{M}_{N-n}} \mathbf{E}_{x}[g(X(\tau
))]$, whence
\[
V_n=\sup_{\tau\in\mathcal{M}_{N-n}} \mathbf{E}_{Z_{n}}[g(X(\tau))],
\]
showing the result.
\end{pf}

\section{Approximation of the value function}
\label{section_Vchap}
To approximate the sequence of value functions $(V_n)$, we proceed in
two steps.
First, the continuous-time maximization of operator $L$ is converted
into a discrete-time maximization by using a path-dependent
time-discretization scheme to give a new operator $L^{d}$. In
particular, it is important to remark that these time-discretization
grids depend on the the post-jump locations $\{Z_{k}\}$ of the PDMP
(see Definition \ref{grid} and Remark \ref{rqgrilletstar}).
Second, the conditional expectations of the Markov chain $(\Theta_{k})$
in the definition of $L^{d}$ are replaced by the conditional
expectations of its quantized approximation $(\widehat{\Theta}_{k})$ to
define an operator $\widehat{L}^{d}$.

First, we define the path-adapted discretization grids as follows.
\begin{definition}
\label{grid}
For $z\in E$, set $\Delta(z)\in\ ]0,t^{*}(z)[$. Define $
n(z)=\operatorname{int} (\frac{t^{*}(z)}{\Delta(z)} )-1$, where
$\operatorname{int}(x)$ denotes the greatest integer smaller than or
equal to $x$.
The set of points $(t_{i})_{i\in\{0,\ldots,n(z)\}}$ with $t_{i} =
i\Delta(z)$ is denoted by $G(z)$.
This is the grid associated with the time interval $[0,t^{*}(z)]$.
\end{definition}
\begin{remark}
\label{rqgrilletstar}
It is important to note that, for all $z\in E$, not only one has
$t^{*}(z) \notin G(z)$, but also $\max G(z)=t_{n(z)}\leq
t^{*}(z)-\Delta
(z)$. This property is crucial for the sequel.
\end{remark}
\begin{definition}
Consider for $w\in\mathbf{B}(E)$ and $z\in E$,
\begin{eqnarray*}
L^{d}(w,g)(z) & = & \max_{s \in G(z)} \bigl\{ \mathbf{E} \bigl[ w(Z_{1})
\mathbf{1}_{\{ S_{1} <s \} }
+g (\phi(z,s) )\mathbf{1}_{\{ S_{1}\geq s\} } | Z_{0}=z
\bigr]
\bigr\} \\
&&{} \vee\mathbf{E} [w(Z_{1}) | Z_{0}=z ] .
\end{eqnarray*}
\end{definition}

Now let us turn to the quantization of $(\Theta_{n})$. The quantization
algorithm will provide us with a finite grid $\Gamma^{\Theta
}_n\subset
E\times\mathbb{R}_+$ at each time $0\leq n\leq N$ as well as weights
for each point of the grid (see, e.g.,
\cite{bally03,pages98,pages04}). Set $p\geq1$ such that $\Theta_n$ has
finite moments at least up to the order $p$ and let $p_n$ be the
closest-neighbor projection from $E\times\mathbb{R}_+$ onto $\Gamma
^{\Theta}_n$ (for the distance of norm $p$; if there are several
equally close neighbors, pick the one with the smallest index). Then
the quantization of $\Theta_n$ is defined by
\[
\widehat{\Theta}_{n}= (\widehat{Z}_{n},\widehat{S}_{n}
)=p_n
({Z}_{n},{S}_{n} ).
\]
We will also denote by $\Gamma^{Z}_{n}$, the projection of $\Gamma
^{\Theta}_n$ on $E$, and by $\Gamma^{S}_{n}$, the projection of
$\Gamma
^{\Theta}_n$ on $\mathbb{R}_+$.

In practice, one will first compute the quantization grids and weights,
and then compute a path-adapted time-grid for each $z\in\Gamma
^{Z}_{n}$, for all $0\leq n\leq N-1$. Hence, there is only a finite
number of time grids to compute, and like the quantization grids, they can
be computed and stored off-line.

The definition of the discretized operators now naturally follows the
characterization given in Lemma \ref{Iproba}.
\begin{definition}
\label{Lchap}
For $k\in\{1,\ldots,N\}$, $w\in\mathbf{B} (\Gamma^{Z}_{k})$, $z\in
\Gamma^{Z}_{k-1}$, and $s\in\mathbb{R}_{+}$
\begin{eqnarray*}
\widehat{J}_{k}(w,g)(z,s) & = & \mathbf{E} \bigl[w(\widehat
{Z}_{k})\mathbf{1}_{\{ \widehat{S}_{k}<s \} }
+g (\phi(z,s) )\mathbf{1}_{ \{ \widehat{S}_{k}\geq s\} }
|
\widehat{Z}_{k-1}=z \bigr], \nonumber\\
\widehat{K}_{k}(w)(z) & = & \mathbf{E} [w(\widehat{Z}_{k}) |\widehat
{Z}_{k-1}=z ], \nonumber\\
\widehat{L}_{k}^{d}(w,g)(z) & = & \max_{s\in G(z)} \{ \widehat
{J}_{k}(w,g)(z,s) \} \vee\widehat{K}_{k}(w)(z).
\end{eqnarray*}
\end{definition}

Note that $\widehat{\Theta}_{n}$ is a random variable taking finitely
many values, hence the expectations above actually are finite sums, the
probability of each atom being given by its weight on the quantization grid.
We can now give the complete construction of the sequence approximating $(V_n)$.
\begin{definition}
\label{vchap}
Consider $\widehat{v}_{N}(z)=g(z)$ where $z\in\Gamma^{Z}_{N}$ and for
$k\in\{1,\break\ldots,N\}$
%
%
\begin{equation}
\widehat{v}_{k-1}(z) = \widehat{L}_{k}^{d}(\widehat{v}_{k},g)(z),
\end{equation}
where $z\in\Gamma^{Z}_{k-1}$.
\end{definition}
\begin{definition}
The approximation of $V_{k}$ is denoted by
\label{Vchap}
%
%
\begin{equation}
\widehat{V}_{k} = \widehat{v}_{k}(\widehat{Z}_{k})
\end{equation}
for $k\in\{0,\ldots,N\}$.
\end{definition}

\section{Error estimation for the value function}
\label{section_erreur}
We are now able to state our main result, namely the convergence of our
approximation scheme with an upper bound for the rate of convergence.
\begin{theorem}
\label{th_Vn}
Set $n\in\{0,\ldots,N-1\}$, and suppose that $\Delta(z)$, for $z\in
\Gamma^{z}_{n}$, are chosen such that
\[
\min_{z\in\Gamma^{z}_{n}}\{\Delta(z)\}>(2C_{\lambda})^{-1/2}
([ t^* ]\|\widehat{Z}_{n}-{Z}_{n}\|_p+\|{S}_{n+1}-\widehat{S}_{n+1}\|
_p )^{1/2}.
\]
Then the discretization error for $V_{n}$ is no greater than the following:
\begin{eqnarray*}
\| V_{n}-\widehat{V}_{n}\|_p
&\leq& \|V_{n+1}-\widehat{V}_{n+1} \|_p+\alpha\|\Delta
(\widehat
{Z}_{n})\|_p+\beta_n\|\widehat{Z}_{n}-{Z}_{n}\|_p\\
&&{}+2[ v_{n+1} ]\|
\widehat{Z}_{n+1}-{Z}_{n+1}\|_p\\
&&{}+\gamma([ t^* ]\|\widehat{Z}_{n}-{Z}_{n}\|_p+\|{S}_{n+1}-\widehat
{S}_{n+1}\|
_p )^{1/2},
\end{eqnarray*}
where $\alpha=[ g ]_{2} +2C_{g} C_{\lambda}$, $\beta_n=[ v_{n} ]+[
v_{n+1} ]_{1}E_2+ C_gE_4+ ( [ g ]_1+[ g ]_2[ t^* ] )\vee
(
[ v_{n+1} ]_*[ Q ] )$, $\gamma=4C_g(2C_{\lambda})^{1/2}$, and
$E_2$ and $E_4$ are defined in Appendix \ref{app_A}.
\end{theorem}

Recall that $V_N=g(Z_N)$ and $\widehat{V}_N=g(\widehat{Z}_{N})$, hence
$\| V_{N}-\widehat{V}_{N}\|_p\leq[g] \|\widehat{Z}_{N}-{Z}_{N} \|_p$.
In addition, the quantization error $\|\Theta_{n}-\widehat{\Theta
}_{n}\|
_p $ goes to zero as the number of points in the grids goes to infinity
(see, e.g.,
\cite{pages98}). Hence $|V_0-\widehat{V}_0|$ can be made arbitrarily
small by
an adequate choice of the discretizations parameters.

Remark that the square root in the last error term is the price to pay
for integrating noncontinuous functions, see the definition of
operator $J$ with the indicator functions, and the introduction of
Section \ref{sectionterme3}.

To prove Theorem \ref{th_Vn}, we split the left-hand side difference
into four terms
\[
\|V_{n}-\widehat{V}_{n}\|_p \leq\sum_{i=1}^{4} \Xi_{i},
\]
where
\begin{eqnarray*}
\Xi_{1} & = & \|v_{n}(Z_{n})-v_{n}(\widehat{Z}_{n})\|_p ,\\
\Xi_{2} & = & \|L(v_{n+1},g)(\widehat
{Z}_{n})-L^{d}(v_{n+1},g)(\widehat
{Z}_{n})\|_p ,\\
\Xi_{3} & = & \|L^{d}(v_{n+1},g)(\widehat{Z}_{n})-\widehat
{L}_{n+1}^{d}(v_{n+1},g)(\widehat{Z}_{n})\|_p ,\\
\Xi_{4} & = & \|\widehat{L}_{n+1}^{d}(v_{n+1},g)(\widehat
{Z}_{n})-\widehat{L}_{n+1}^{d}(\widehat{v}_{n+1},g)(\widehat
{Z}_{n})\|_p.
\end{eqnarray*}
The first term is easy enough to handle thanks to Proposition \ref
{prop1} in Appendix \ref{appvnlip}.
\begin{lemma}
\label{lemerreur1}
A upper bound for $\Xi_{1}$ is
\[
\|v_{n}(Z_{n})-v_{n}(\widehat{Z}_{n})\|_p\leq[v_{n}] \|Z_{n}-\widehat
{Z}_{n}\|_p.
\]
\end{lemma}

We are going to study the other terms one by one in the following sections.

\subsection{Second term}
In this part we study the error induced by the replacement of the
supremum over all nonnegative $t$ smaller than or equal to $t^{*}(z)$
by the maximum over the finite grid $G(z)$ in the definition of
operator $L$.\vadjust{\goodbreak}
\begin{lemma}
Let $w\in\mathbf{L}^{c}$. Then for all $z\in E$,
\[
\Bigl|\sup_{t\leq t^*(z)}J(w,g)(z,t)-\max_{s\in G(z)}J(w,g)(z,s) \Bigr|
\leq( C_wC_{\lambda}+[ g ]_2+C_{g}C_{\lambda} ) \Delta(z).
\]
\end{lemma}
\begin{pf}
Clearly, there exists $\overline{t} \in[0,t^{*}(z)]$ such that $\sup
_{t\leq t^*(z)}J(w,g)(z,t)=J(w,g)(z,\overline{t})$, and there exists
$0\leq i\leq n(z)$ such that $\overline{t} \in[t_i,t_{i+1}]$ [with
$t_{n(z)+1}=t^{*}(z)$]. Consequently, Lemma \ref{Tech2} yields
\begin{eqnarray*}
0 &\leq& \sup_{t\leq t^*(z)}J(w,g)(z,t)-\max_{s\in
G(z)}J(w,g)(z,s)\\
&\leq&
J(w,g)(z,\overline{t})-J(w,g)(z,t_i)\\
&\leq& ( C_{w} C_{\lambda} + [ g ]_{2} +C_{g} C_{\lambda} )
|\overline{t}-t_i|\\
&\leq& ( C_{w} C_{\lambda} + [ g ]_{2} +C_{g} C_{\lambda} )
|t_{i+1}-t_i|,
\end{eqnarray*}
implying the result.
\end{pf}

Turning back to the second error term, one gets the following bound.
\begin{lemma}
\label{lemerreur2}
A upper bound for $\Xi_{2}$ is
\[
\|L(v_{n+1},g)(\widehat{Z}_{n})-L^{d}(v_{n+1},g)(\widehat{Z}_{n})\|_p
\leq( [ g ]_{2} +2C_{g} C_{\lambda} )\|\Delta(\widehat
{Z}_{n})\|_p.
\]
\end{lemma}
\begin{pf}
From the definition of $L$ and $L^{d}$ we readily obtain
\begin{eqnarray*}
&&
\|L(v_{n+1},g) (\widehat{Z}_{n})-L^{d}(v_{n+1},g)(\widehat{Z}_{n})\|
_p \\
&&\qquad\leq\Bigl\|\sup_{t\leq t^*(\widehat{Z}_{n})}J(v_{n+1},g)(\widehat
{Z}_{n},t)-\max_{s\in G(\widehat{Z}_{n})}J(v_{n+1},g)(\widehat
{Z}_{n},s) \Bigr\|_p.
\end{eqnarray*}
Now in view of the previous lemma, one has
\begin{eqnarray*}
&&
\|L(v_{n+1},g) (\widehat{Z}_{n})-L^{d}(v_{n+1},g)(
\widehat{Z}_{n})\|
_p \\
&&\qquad\leq( C_{v_{n+1}} C_{\lambda} + [ g ]_{2} +C_{g} C_{\lambda}
)\|\Delta(\widehat{Z}_{n})\|_p.
\end{eqnarray*}
Finaly, note that $C_{v_{n+1}}= C_{g}$ (see Appendix \ref{appvnlip}),
completing the proof.
\end{pf}

\subsection{Third term}
\label{sectionterme3}
This is the crucial part of our derivation, where we need to compare
conditional expectations relative to the real Markov chain $(Z_n,S_n)$
and its quantized approximation $(\widehat{Z}_n,\widehat{S}_n)$. The
main difficulty stems from the fact that some functions inside the
expectations are indicator functions and in particular they are not
Lipschitz-continuous. We manage to overcome this difficulty by proving
that the event on which the discontinuity actually occurs is of small
enough probability; this is the aim of the following two lemmas.
\begin{lemma}
\label{lem_indic}
For all $n\in\{0,\ldots,N-1\}$ and $0<\eta<\min_{z\in\Gamma
^Z_{n}}\{\Delta(z)\}$,
\begin{eqnarray*}
&&
\Bigl\| \max_{s\in G(\widehat{Z}_{n})}\mathbf{E} \bigl[\bigl|\mathbf
{1}_{\{
{S}_{n+1}<s\}}-\mathbf{1}_{\{\widehat{S}_{n+1}<s\}}\bigr| | \widehat{Z}_{n}
\bigr]
\Bigr\|_p\\
&&\qquad \leq \frac{2}{\eta}\|{S}_{n+1}-\widehat{S}_{n+1}\|_p +
C_{\lambda
}\eta+
\frac{2[t^*] \|Z_{n}-\widehat{Z}_{n}\|_p}{\eta}.
\end{eqnarray*}
\end{lemma}
\begin{pf}
Set $0<\eta<\min_{z\in\Gamma^Z_{n}}\{\Delta(z)\}$. Remark that the
difference of indicator functions is nonzero if and only if
${S}_{n+1}$ and $\widehat{S}_{n+1}$ are on either side of $s$. Hence,
one has
\[
\bigl|\mathbf{1}_{\{{S}_{n+1}<s\}}-\mathbf{1}_{\{\widehat{S}_{n+1}<s\}
}
\bigr|\leq\mathbf{1}_{\{|{S}_{n+1}-\widehat{S}_{n+1}|>{\eta}/{2}\}} +
\mathbf
{1}_{\{|{S}_{n+1}-s|\leq{\eta}/{2}\}}.
\]
This yields
%
%
\begin{eqnarray}\label{lem_indic1}\quad
&& \Bigl\|\max_{s\in G(\widehat{Z}_{n})}\mathbf{E}
\bigl[\bigl|\mathbf
{1}_{\{{S}_{n+1}<s\}}-\mathbf{1}_{\{\widehat{S}_{n+1}<s\}}\bigr| |
\widehat{Z}_{n} \bigr] \Bigr\|_p\nonumber\\[-8pt]\\[-8pt]
&&\qquad\leq \bigl\|\mathbf{1}_{\{|{S}_{n+1}-\widehat{S}_{n+1}|>
{\eta}/{2}\}}
\bigr\|_p+ \Bigl\|\max_{s\in G(\widehat{Z}_{n})}
\mathbf{E} \bigl[\mathbf{1}_{\{s-{\eta}/{2}\leq{S}_{n+1}\leq
s+{\eta}/{2}\}} | \widehat{Z}_{n} \bigr] \Bigr\|_p.\nonumber
\end{eqnarray}
On the one hand, Chebyshev's inequality yields
%
%
\begin{equation}\label{lem_indic2}\qquad
\bigl\| \mathbf{1}_{\{|{S}_{n+1}-\widehat{S}_{n+1}|>{\eta}/{2}\}}
\bigr\|_{p}^{p} =
\mathbf{P}\biggl( |{S}_{n+1}-\widehat{S}_{n+1}|>\frac{\eta}{2} \biggr) \leq
\frac{2^{p}
\| {S}_{n+1}-\widehat{S}_{n+1} \|_{p}^{p}}{\eta^p}.
\end{equation}
On the other hand, as $s\in G(\widehat{Z}_{n})$ and by definition of
$\eta$, one has $s+\eta<t^{*}(\widehat{Z}_{n})$ (see Remark \ref
{rqgrilletstar}). Thus, one has
%
%
\begin{eqnarray}\label{lem_indic3}
&&
\mathbf{E} \bigl[\mathbf{1}_{\{s-{\eta}/{2}\leq{S}_{n+1}\leq
s+{\eta}/{2}\}} | \widehat{Z}_{n} \bigr]
\nonumber\\
&&\qquad=\mathbf{E} \bigl[\mathbf{E} \bigl[\mathbf{1}_{\{s-{\eta
}/{2}\leq
{S}_{n+1}\leq s+{\eta}/{2}\}} | Z_{n} \bigr] | \widehat{Z}_{n}
\bigr]\nonumber\\[-8pt]\\[-8pt]
&&\qquad\leq\mathbf{E} \biggl[\int_{s-{\eta}/{2}}^{s+{\eta}/{2}}\lambda
(\phi(Z_{n},u))\,du | \widehat{Z}_{n} \biggr]
+ \mathbf{E} \bigl[ \mathbf{1}_{\{ t^{*}(Z_{n})\leq s+{\eta}/{2} \}
} | \widehat{Z}_{n} \bigr] \nonumber\\
&&\qquad\leq\eta C_{\lambda} + \mathbf{E} \bigl[ \mathbf{1}_{\{
t^{*}(Z_{n})\leq t^{*}(\widehat{Z}_{n})-{\eta}/{2} \}} |
\widehat{Z}_{n}\bigr].\nonumber
\end{eqnarray}
Combining equations (\ref{lem_indic1})--(\ref{lem_indic3}), the result
follows.
\end{pf}
\begin{lemma}
\label{lem_indic4}
For all $n\in\{0,\ldots,N\}$ and $0<\eta<\min_{z\in\Gamma
^Z_{n}}\{
\Delta(z)\}$,
\[
\bigl\|\mathbf{1}_{t^{*}(Z_{n})<t^{*}(\widehat{Z}_{n})-\eta} \bigr\|_p
\leq
\frac{[t^*] \|Z_{n}-\widehat{Z}_{n}\|_p}{\eta}.
\]
\end{lemma}
\begin{pf}
We use Chebyshev's inequality again. One clearly has
\begin{eqnarray*}
\mathbf{E} \bigl[\bigl|\mathbf{1}_{t^{*}(Z_{n})<t^{*}(\widehat{Z}_{n})-\eta
}\bigr|^p
\bigr] &=&
\mathbf{P} \bigl( t^{*}(Z_{n})<t^{*}(\widehat{Z}_{n})-\eta\bigr) \\
&\leq&\mathbf{P} \bigl( | t^{*}(Z_k)-t^{*}(\widehat{Z}_k)
|>\eta
\bigr) \\
&\leq&\frac{[t^*]^p\|Z_k-\widehat{Z}_k\|_p^p}{\eta^p},
\end{eqnarray*}
showing the result.
\end{pf}

Now we turn to the consequences of replacing the Markov chain
$(Z_n,S_n)$ by its quantized approximation $(\widehat{Z}_n,\widehat
{S}_n)$ in the conditional expectations.
\begin{lemma}\label{lemKquantif}
Let $w\in\mathbf{L}^{c}$, then one has
\begin{eqnarray*}
&&
| \mathbf{E} [w({Z}_{n+1}) |{Z}_{n}=\widehat{Z}_{n}
]-\mathbf
{E} [w(\widehat{Z}_{n+1}) |\widehat{Z}_{n} ] | \\[-0.5pt]
&&\qquad\leq( C_{w} E_{4} + [ w ]_{1} E_{2} + [ w ]_{*} [ Q ] )
\mathbf{E} [|{Z}_{n}-\widehat{Z}_{n}| | \widehat{Z}_{n} ]\\[-0.5pt]
&&\qquad\quad{} +[ w ]\mathbf
{E} [ | Z_{n+1}-\widehat{Z}_{n+1} | | \widehat{Z}_{n} ].
\end{eqnarray*}
\end{lemma}
\begin{pf}
First, note that\vspace*{-1pt}
\begin{eqnarray*}
&&
\mathbf{E} [w({Z}_{n+1}) |{Z}_{n}=\widehat{Z}_{n} ]-\mathbf
{E}
[w(\widehat{Z}_{n+1}) |\widehat{Z}_{n} ]
\\[-0.5pt]
&&\qquad= \mathbf{E} [w({Z}_{n+1}) |{Z}_{n}=\widehat{Z}_{n}
]-\mathbf
{E} [w({Z}_{n+1}) |\widehat{Z}_{n} ]\\[-0.5pt]
&&\qquad\quad{} + \mathbf{E} [w({Z}_{n+1}) |\widehat{Z}_{n} ]-\mathbf{E}
[w(\widehat{Z}_{n+1}) |\widehat{Z}_{n} ].
\end{eqnarray*}
On the one hand, Remark \ref{remarque_lip} yields\vspace*{-0.5pt}
\[
| \mathbf{E} [w({Z}_{n+1}) | \widehat{Z}_{n} ] - \mathbf
{E}
[w(\widehat{Z}_{n+1}) | \widehat{Z}_{n} ] |
\leq[ w ]\mathbf{E} [ | Z_{n+1}-\widehat{Z}_{n+1} | | \widehat
{Z}_{n} ].
\]
On the other hand, recall that by construction of the quantized
process, one has $ (\widehat{Z}_{n},\widehat{S}_{n} )=p_n
({Z}_{n},{S}_{n} )$.
Hence we have the following property: $\sigma\{\widehat{Z}_{n}\}
\subset
\sigma\{{Z}_{n}, S_n\}$.
By using the special structure of the PDMP $\{X(t)\}$, we have
$\sigma\{{Z}_{n}, S_n\} \subset\mathcal{F}_{T_{n}}$.
Now, by using the Markov property of the process $\{X(t)\}$, it follows
that\vspace*{-1pt}
\[
\mathbf{E} [w({Z}_{n+1}) |\widehat{Z}_{n} ] = \mathbf{E} [
\mathbf{E} [w({Z}_{n+1}) | \mathcal{F}_{T_{n}} ] |
\widehat{Z}_{n} ] = \mathbf{E} [ \mathbf{E} [w({Z}_{n+1}) | Z_{n}
] | \widehat{Z}_{n} ] .
\]
Equation (\ref{Kproba}) thus yields\vspace*{-0.5pt}
\begin{eqnarray*}
&&\mathbf{E} [w({Z}_{n+1}) |{Z}_{n}=\widehat{Z}_{n}
]-\mathbf
{E} [w({Z}_{n+1}) |\widehat{Z}_{n} ]\\[-0.5pt]
&&\qquad=\mathbf{E} \bigl[\mathbf{E} [w({Z}_{n+1}) |{Z}_{n}=\widehat{Z}_{n}
]-\mathbf{E} [w({Z}_{n+1}) |{Z}_{n} ] |\widehat{Z}_{n}
\bigr]\\[-0.5pt]
&&\qquad=\mathbf{E} [Kw(\widehat{Z}_{n})-Kw({Z}_{n}) |\widehat{Z}_{n} ].
\end{eqnarray*}
Now we use Lemma \ref{Tech4} to conclude.
\end{pf}

Now we combine the preceding lemmas to derive the third error term.
\begin{lemma}
\label{lemerreur3}
For all $0<\eta<\min_{z\in\Gamma^Z_{n}}\{\Delta(z)\}$, an upper
bound for $\Xi_{3}$ is
\begin{eqnarray*}
&&
\| L^{d} (v_{n+1},g)(\widehat{Z}_{n})-\widehat
{L}_{n+1}^{d}(v_{n+1},g)(\widehat{Z}_{n}) \|_p\\[-1pt]
&&\qquad
\leq\biggl\{[ v_{n+1} ]_{1}E_2+ C_gE_4+ 2C_g\frac{[ t^* ]}{\eta
}\\[-2.5pt]
&&\qquad\quad\hspace*{4.1pt}{}+
( [ g ]_1+[ g ]_2[ t^* ] )\vee( [ v_{n+1} ]_*[ Q ] )
\biggr\}\|\widehat{Z}_{n}-{Z}_{n}\|_p\\[-2.5pt]
&&\qquad\quad{} +[ v_{n+1} ]\|\widehat{Z}_{n+1}-{Z}_{n+1}\|_p+2C_g \biggl(2C_{\lambda
}\eta+\frac{\|{S}_{n+1}-\widehat{S}_{n+1}\|_p}{\eta} \biggr).\
\end{eqnarray*}
\end{lemma}
\begin{pf}
To simplify notation, set $\Psi
(x,y,t)=v_{n+1}(y)\mathbf{1}_{\{t<s\}}+g (\phi(x,t) )\times\break
\mathbf{1}_{\{t\geq s\}}$.
From the definition of $L^{d}$ and $\widehat{L}_{n+1}^{d}$, one
readily obtains
%
%
\begin{eqnarray}\label{Ld}
&&| L^{d}(v_{n+1},g)(\widehat{Z}_{n})-\widehat
{L}_{n+1}^{d}(v_{n+1},g)(\widehat{Z}_{n}) |
\nonumber\\
&&\qquad\leq\max_{s\in G(\widehat{Z}_{n})} |\mathbf{E} [\Psi(Z_n,
Z_{n+1}, S_{n+1}) | {Z}_{n}=\widehat{Z}_{n} ]\nonumber\\[-8pt]\\[-8pt]
&&\qquad\quad\hspace*{50pt}{}-\mathbf{E}
[\Psi(\widehat{Z}_n, \widehat{Z}_{n+1}, \widehat{S}_{n+1}) |
\widehat{Z}_{n}]
|\nonumber\\
&&\qquad\quad{}\vee|\mathbf{E} [v_{n+1}({Z}_{n+1})
|{Z}_{n}=\widehat{Z}_{n} ]-\mathbf{E} [v_{n+1}(\widehat{Z}_{n+1})
|\widehat{Z}_{n} ] |.\nonumber
\end{eqnarray}
On the one hand, combining Lemma \ref{lemKquantif} and the fact that
$v_{n+1}$ is in $\mathbf{L}^{c}$ 
(see Proposition \ref{prop1}), we obtain
%
%
\begin{eqnarray}\label{Ldchap2solved}
&&
|\mathbf{E} [v_{n+1}({Z}_{n+1}) |{Z}_{n}=\widehat{Z}_{n}
]-\mathbf{E} [v_{n+1}(\widehat{Z}_{n+1}) |\widehat{Z}_{n} ]
|\nonumber\\
&&\qquad\leq [ v_{n+1} ] \mathbf{E} [|{Z}_{n+1}-\widehat{Z}_{n+1} |
\widehat{Z}_{n} ]\\
&&\qquad\quad{} + ( C_{g} E_{4} +[ v_{n+1} ]_{1} E_{2}+[ v_{n+1} ]_{*}[ Q ]
) \mathbf{E} [|{Z}_{n}-\widehat{Z}_{n}| | \widehat{Z}_{n}].\nonumber
\end{eqnarray}
On the other hand, similar arguments as in the proof of
Lemma \ref{lemKquantif} yield
%
%
\begin{eqnarray}\label{Ld1b}
&&
\mathbf{E} [\Psi(Z_n, Z_{n+1}, S_{n+1}) |
{Z}_{n}=\widehat{Z}_{n} ]-\mathbf{E} [\Psi(\widehat{Z}_n, \widehat
{Z}_{n+1}, \widehat{S}_{n+1}) | \widehat{Z}_{n} ]\nonumber\\
&&\qquad=\mathbf{E} \bigl[\mathbf{E} [\Psi(Z_n, Z_{n+1}, S_{n+1}) |
{Z}_{n}=\widehat{Z}_{n} ]\nonumber\\
&&\qquad\quad\hspace*{8.6pt}{}-\mathbf{E} [\Psi(Z_n, Z_{n+1},
S_{n+1}) | {Z}_{n}=Z_n ] |\widehat{Z}_{n} \bigr]\\
&&\qquad\quad{}+\mathbf{E} [\Psi(Z_n, Z_{n+1}, S_{n+1}) | \widehat{Z}_{n}
]-\mathbf{E} [\Psi(\widehat{Z}_n, \widehat{Z}_{n+1}, \widehat
{S}_{n+1}) | \widehat{Z}_{n} ]\nonumber\\
&&\qquad=\Upsilon_{1}+\Upsilon_{2} . \nonumber
\end{eqnarray}
The second difference of the right-hand side of (\ref{Ld1b}),
labeled $\Upsilon_{2}$, clearly satisfies
%
%
\begin{eqnarray}\label{Ld1bsolved}
|\Upsilon_{2}| & \leq& [ v_{n+1} ] \mathbf{E} [|\widehat
{Z}_{n+1}-{Z}_{n+1}| | \widehat{Z}_{n} ]+[ g ]_1\mathbf{E}
[|\widehat{Z}_{n}-{Z}_{n}| | \widehat{Z}_{n} ]\nonumber\\[-8pt]\\[-8pt]
&&{} + 2C_g\mathbf{E} \bigl[\bigl|\mathbf{1}_{\{{S}_{n+1}<s\}}-\mathbf{1}_{\{
\widehat{S}_{n+1}<s\}}\bigr| | \widehat{Z}_{n} \bigr].\nonumber
\end{eqnarray}
Let us turn now to the first difference of the right-hand side of
(\ref{Ld1b}), labeled $\Upsilon_{1}$. We meet another
difficulty here. Indeed, we know by construction that $s<t^{*}(\widehat
{Z}_n)$, but we know nothing regarding the relative positions of $s$
and $t^{*}({Z}_n)$. In the event where $s\leq t^{*}({Z}_n)$ as well,
we recognize operator $J$ inside the expectations. In the opposite
event $s>t^{*}({Z}_n)$, we crudely bound $\Psi$ by
$C_{v_{n+1}}+C_g=2C_g$. Hence, one obtains
\begin{eqnarray*}
|\Upsilon_{1}| &\leq& \mathbf{E} \bigl[ | J(v_{n+1},g)(\widehat
{Z}_{n},s)-J(v_{n+1},g)({Z}_{n},s) |
\mathbf{1}_{\{ s\leq t^{*}({Z}_{n}) \} } |\widehat{Z}_{n} \bigr]
\\
&&{}+2C_g\mathbf{E} \bigl[\mathbf{1}_{ \{ t^{*}(Z_{n})<s \} } | \widehat
{Z}_{n} \bigr].
\end{eqnarray*}
Now Lemma \ref{Tech3} gives an upper bound for the first term. As for
the indicator function, by definition of $G(\widehat{Z}_{n})$ and our choice
of $\eta$, we have $s<t^{*}(\widehat{Z}_{n})-\eta$. Thus, one has
%
%
\begin{eqnarray}\label{Ld1asolved}
|\Upsilon_{1}| &\leq& (C_{g}E_{1}+[ v_{n+1} ]_{1}E_{2}+E_{3}
) \mathbf{E} [|\widehat{Z}_{n}-{Z}_{n}| |
\widehat{Z}_{n}]\nonumber\\[-8pt]\\[-8pt]
&&{} + 2C_g\mathbf
{E} \bigl[\mathbf{1}_{ \{ t^{*}(Z_{n})<t^{*}(\widehat{Z}_{n})-\eta\} }
|\widehat{Z}_{n} \bigr].\nonumber
\end{eqnarray}
Now, combining (\ref{Ld}), (\ref{Ldchap2solved}), (\ref{Ld1bsolved})
and (\ref{Ld1asolved}), and the fact that $C_g E_1+E_3= C_g E_4+[ g
]_1+[ g ]_2[ t^* ]$, one gets
\begin{eqnarray*}
&&
| L^{d}(v_{n+1}, g) (\widehat{Z}_{n})-\widehat
{L}_{n+1}^{d}(v_{n+1},g)(\widehat{Z}_{n}) | \\
&&\qquad\leq\{[ v_{n+1} ]_{1}E_2+ C_gE_4\\
&&\qquad\quad\hspace*{1.7pt}{} + ( [ g ]_1+[ g ]_2[ t^* ] )\vee( [
v_{n+1} ]_*[ Q ] ) \}\mathbf
{E} [|\widehat{Z}_{n}-{Z}_{n}| | \widehat{Z}_{n} ]\\
&&\qquad\quad{} +[ v_{n+1} ]\mathbf{E} [|\widehat{Z}_{n+1}-{Z}_{n+1}| |
\widehat{Z}_{n} ] \\
&&\qquad\quad{} +2C_g\mathbf{E} \bigl[\mathbf{1}_{t^{*}(Z_{n})<t^{*}(\widehat
{Z}_{n})-\eta} | \widehat{Z}_{n} \bigr]\\
&&\qquad\quad{} +2C_g\max_{s\in G(\widehat
{Z}_{n})}\mathbf{E} \bigl[\bigl|\mathbf{1}_{\{{S}_{n+1}<s\}}-\mathbf{1}_{\{
\widehat{S}_{n+1}<s\}}\bigr| | \widehat{Z}_{n} \bigr].
\end{eqnarray*}
Finally, we conclude by taking the $L^p$ norm on both sides and using
Lemmas \ref{lem_indic} and \ref{lem_indic4}.
\end{pf}

\subsection{Fourth term}
The last error term is a mere comparison of two finite sums.
\begin{lemma}
\label{lemerreur4}
An upper bound for $\Xi_{4}$ is
\begin{eqnarray*}
&&
\|\widehat{L}_{n+1}^{d}(v_{n+1},g)(\widehat{Z}_{n})-\widehat
{L}_{n+1}^{d}(\widehat{v}_{n+1},g)(\widehat{Z}_{n})\|_p\\
&&\qquad\leq[ v_{n+1}
] \|\widehat{Z}_{n+1}-{Z}_{n+1} \|_p + \|
V_{n+1}-\widehat
{V}_{n+1} \|_p.
\end{eqnarray*}
\end{lemma}
\begin{pf}
By definition of operator $\widehat
{L}_{n}^{d}$, one has
\begin{eqnarray*}
&&
\|\widehat{L}_{n+1}^{d}(v_{n+1},g)(\widehat
{Z}_{n})-\widehat
{L}_{n+1}^{d}(\widehat{v}_{n+1},g)(\widehat{Z}_{n})\|_p\\
&&\qquad= \Bigl\|\max_{s\in G(\widehat{Z}_{n})} \bigl\{ \mathbf{E}
\bigl[v_{n+1}(\widehat{Z}_{n+1})\mathbf{1}_{\{ \widehat{S}_{n+1}<s \}
}+g(\phi(\widehat{Z}_{n},s) )\mathbf{1}_{ \{ \widehat{S}_{n+1}\geq
s\}} | \widehat{Z}_{n} \bigr] \bigr\}\\
&&\qquad\quad\hspace*{2.1pt}{}\vee\mathbf{E}
[v_{n+1}(\widehat{Z}_{n+1}) |\widehat{Z}_{n} ]\\
&&\qquad\quad\hspace*{2.1pt}{}-\max_{s\in G(\widehat{Z}_{n})} \bigl\{ \mathbf{E} \bigl[\widehat
{v}_{n+1}(\widehat{Z}_{n+1})\mathbf{1}_{\{ \widehat{S}_{n+1}<s \}
}\\
&&\qquad\quad\hspace*{61.3pt}{}+g (\phi(\widehat{Z}_{n},s) )\mathbf{1}_{ \{ \widehat
{S}_{n+1}\geq s\} } | \widehat{Z}_{n} \bigr] \bigr\}\vee
\mathbf{E} [\widehat{v}_{n+1}(\widehat{Z}_{n+1}) |\widehat{Z}_{n} ] \Bigr\|
_p\\
&&\qquad\leq \|\mathbf{E}[v_{n+1}(\widehat{Z}_{n+1})-\widehat
{v}_{n+1}(\widehat{Z}_{n+1}) |\widehat{Z}_{n}] \|_p\\
&&\qquad\leq \|v_{n+1}(\widehat{Z}_{n+1})-v_{n+1}({Z}_{n+1}) \|_p+ \|
v_{n+1}({Z}_{n+1})-\widehat{v}_{n+1}(\widehat{Z}_{n+1}) \|_p.
\end{eqnarray*}
We conclude using the fact that $v_{n+1} \in\mathbf{L}^{c}$ (see
Proposition \ref{prop1}) and the definitions of $V_{n}$ and $\widehat{V}_{n}$.
\end{pf}

\subsection[Proof of Theorem 5.1]{Proof of Theorem \protect\ref{th_Vn}}
We can finally turn to the proof of Theorem \ref{th_Vn}. Lemmas \ref
{lemerreur1}, \ref{lemerreur2}, \ref{lemerreur3} and \ref{lemerreur4}
from the preceding sections directly yield, for all $
0<\eta<\min_{z\in\Gamma^{z}_{n}}\{\Delta(z)\}$,
\begin{eqnarray*}
\|V_{n}-\widehat{V}_{n}\|_p
&\leq&
[ v_{n} ] \|\widehat{Z}_{n}-{Z}_{n}\|_p+ ([ g ]_{2} +2C_{g}
C_{\lambda} )\|\Delta(\widehat{Z}_{n})\|_p\\
&&{} + \biggl\{[ v_{n+1} ]_{1}E_2+ C_gE_4+ 2C_g\frac{[ t^* ]}{\eta
}\\
&&\hspace*{17.5pt}{}+ ( [ g ]_1+[ g ]_2[ t^* ] )\vee( [ v_{n+1} ]_*[ Q ]
) \biggr\}\|\widehat{Z}_{n}-{Z}_{n}\|_p\\
&&{} +[ v_{n+1} ]\|\widehat{Z}_{n+1}-{Z}_{n+1}\|_p+2C_g \biggl(2C_{\lambda}\eta
+\frac{\|{S}_{n+1}-\widehat{S}_{n+1}\|_p}{\eta} \biggr)\\
&&{} +[ v_{n+1} ] \|\widehat{Z}_{n+1}-{Z}_{n+1} \|_p+ \|
V_{n+1}-\widehat{V}_{n+1} \|_p.
\end{eqnarray*}
The optimal choice for $\eta$ clearly satisfies
\[
2C_{\lambda}\eta=\frac{1}{\eta} ([ t^* ]\|\widehat{Z}_{n}-{Z}_{n}\|
_p+\|
{S}_{n+1}-\widehat{S}_{n+1}\|_p ),
\]
providing it also satisfies the condition $0<\eta<\min_{z\in
\Gamma
^{z}_{n}}\{\Delta(z)\}$. Hence, rearranging the terms above, one gets
the expected result
\begin{eqnarray*}
\| V_{n}-\widehat{V}_{n}\|_p &\leq& \|V_{n+1}-\widehat
{V}_{n+1}
\|_p+ ([ g ]_{2} +2C_{g} C_{\lambda} ) \|\Delta(\widehat
{Z}_{n})\|_p\\
&&{}+ \{[ v_{n} ]+[ v_{n+1} ]_{1}E_2+ C_gE_4\\
&&\hspace*{15.4pt}{}+ ( [ g ]_1+[ g ]_2[ t^* ]
)\vee( [ v_{n+1} ]_*[ Q ] ) \}\|\widehat{Z}_{n}-{Z}_{n}\|_p\\
&&{}+2[ v_{n+1} ]\|\widehat{Z}_{n+1}-{Z}_{n+1}\|_p\\
&&{}+4C_g(2C_{\lambda
})^{1/2}
([ t^* ]\|\widehat{Z}_{n}-{Z}_{n}\|_p+\|{S}_{n+1}-\widehat{S}_{n+1}\|
_p )^{1/2}.
\end{eqnarray*}

\section{Numerical construction of an $\epsilon$-optimal stopping time}
\label{sectiontpsarret}
In \cite{gugerli86}, Theorem 1, Gugerli defined an $\epsilon$-optimal
stopping time for the original problem.
Roughly speaking, this stopping time depends on the embedded Markov
chain $(\Theta_n)$ and on the optimal value function.
Therefore, a natural candidate for an $\epsilon$-optimal stopping time
should be obtained by replacing the Markov chain $(\Theta_n)$ and the
optimal value function by their \textit{quantized approximations}.
However, this leads to un-tractable comparisons between some quantities
involving $(\Theta_n)$ and its quantized approximation.
It is then far from obvious to show that this method would provide a
computable $\epsilon$-optimal stopping rule.
Nonetheless, by modifying the approach of Gugerli \cite{gugerli86}, we
are able to propose a numerical construction of an $\epsilon$-optimal
stopping time of the original stopping problem.

Here is how we proceed. First, recall that $p_n$ be the
closest-neighbor projection from $E\times\mathbb{R}_+$ onto $\Gamma
^{\Theta}_n$, and for all $(z,s)\in E\times\mathbb{R}_+$ define
$(\widehat{z}_n,\widehat{s}_n)=p_n(z,s)$. Note that $\widehat{z}_n$
and $\widehat{s}_n$ depend on
both $z$ and $s$. Now, for $n\in\{1,\ldots,N\}$, define
\[
s^{*}_{n}(z,s) = \min\Bigl\{ t\in G(\widehat{z}_{n-1}) | \widehat
{J}_{n}(\widehat{v}_{n},g)(\widehat{z}_{n-1},t)=\max_{u\in
G(\widehat{z}_{n-1})}
\widehat{J}_{n}(\widehat{v}_{n},g)(\widehat{z}_{n-1},u) \Bigr\}
\]
and
\[
r_{n,\beta}(z,s)=
\cases{
t^{*}(z),\qquad\mbox{if $\displaystyle\widehat{K}_{n}\widehat{v}_{n}(\widehat{z}_{n-1})>
\max
_{u\in G(\widehat{z}_{n-1})} \widehat{J}_{n}(\widehat
{v}_{n},g)(\widehat{z}_{n-1},u)$},
\cr
s^{*}_{n}(z,s) \mathbf{1}_{\{s^{*}_{n}(z,s)<t^{*}(z)\}} +
\bigl(t^{*}(z)-\beta\bigr) \mathbf{1}_{\{s^{*}(z,s)\geq t^{*}(z)\}},\cr
\hspace*{49pt} \mbox{otherwise.}}
\]
Note the use of both the real jump time horizon $t^{*}(z)$ and the
quantized approximations of $K$, $J$ and $(z,s)$.
Set
\[
\tau_{1} = r_{N,\beta}(Z_{0}, S_0) \wedge T_{1}
\]
and for $n\in\{1,\ldots,N-1\}$, set
\[
\tau_{n+1} =
\cases{
r_{N-n, \beta}(Z_{0},S_0), &\quad if $T_{1}> r_{N-n,\beta}(Z_{0},S_0)$,
\cr
T_{1}+\tau_{n}\circ\theta_{T_{1}}, &\quad otherwise.}
\]

Our stopping rule is then defined by $\tau_{N}$.
\begin{remark}
This procedure is especially appealing because it requires no more
calculation: we have already computed the values of $\widehat{K}_{n}$
and $\widehat{J}_{n}$ on the grids. One just has to store the point
where the maximum of $\widehat{J}_{n}$ is reached.
\end{remark}
\begin{lemma}
$\tau_{N}$ is an $\{\mathcal{F}_{T}\}$-stopping time.
\end{lemma}
\begin{pf}
Set $U_{1}=r_{1,\beta}(Z_{0},S_{0})$ and for $2\leq k \leq N$
$U_{k}=r_{k,\beta}(Z_{k-1},S_{k-1}) \times\break\mathbf{1}_{\{r_{k-1,\beta
}(Z_{k-2},S_{k-2})\geq S_{k-1}\}}$.
One then clearly has $\tau_{N}=\sum_{k=1}^{N}U_{k}\wedge S_{k}$
which is an $\{\mathcal{F}_{T}\}$-stopping time by Proposition
\ref{Aprop2}.
\end{pf}

Now let us show that this stopping time provides a good approximation
of the value function $V_{0}$.
Namely, for all $z\in E$ set
\[
\overline{v}_{n}(z) = \mathbf{E} [ g(X_{\tau_{N-n}}) |
Z_{n}=z ]
\]
and in accordance to our previous notation introduce, for $n\in\{
1,\ldots,N-1\}$
\[
\overline{V}_{n} = \overline{v}_{n}(Z_{n}).
\]
The comparison between $V_{0}$ and $\overline{V}_{0}$ is provided by the
next two theorems.
\begin{theorem}\label{thtpsarret}
Set $n\in\{0,\ldots,N-2\}$ and suppose the discretization parameters
are chosen such that there exists $0<a<1$ satisfying
\[
\frac{\beta}{a}=(2C_{\lambda})^{-1/2} \biggl(\frac{[ t^* ]}{1-a}\|
\widehat{Z}_{n}-{Z}_{n}\|_p+\|{S}_{n+1}-\widehat{S}_{n+1}\|_p \biggr)^{1/2}
<\min
_{z\in
\Gamma^{z}_{n}}\{\Delta(z)\}.
\]
Then one has
\begin{eqnarray*}
\| \overline{V}_{n} - V_{n} \|_{p} & \leq& \| \overline{V}_{n+1} -
V_{n+1} \|_{p} + \| \widehat{V}_{n+1} -
V_{n+1} \|
_{p}+ \| \widehat{V}_{n} - V_{n} \|_{p}\\
&&{}+2[ v_{n+1} ] \| {Z}_{n+1}-\widehat{Z}_{n+1} \|_p +a_n \|
{Z}_{n}-\widehat{Z}_{n} \|_p\\
&&{}+4C_g(2C_{\lambda})^{1/2} \biggl(\frac{[ t^* ]}{1-a}\|\widehat
{Z}_{n}-{Z}_{n}\|_p+\|{S}_{n+1}-\widehat{S}_{n+1}\|_p \biggr)^{1/2}
\end{eqnarray*}
with $a_n= (2[ v_{n+1} ]_1E_2+2C_gC_{t^*}[ \lambda
]_1(2+C_{t^*}C_{\lambda})+ (4C_gC_{\lambda}[ t^* ]+2[ v_{n+1} ]_*[ Q
] )\vee(3[ g ]_1) )$.
\end{theorem}
\begin{pf}
The definition of $\tau_{n}$ and the strong Markov property of the
process $\{X(t)\}$ yield
\begin{eqnarray*}
\overline{v}_{n}(Z_{n})
&= & \mathbf{E} \bigl[ g\bigl(X_{r_{n+1,\beta}(Z_{n},S_n)}\bigr) \mathbf{1}_{\{
S_{n+1}> r_{n+1,\beta}(Z_{n},S_n)\}} | Z_{n} \bigr]\\
&&{}+\mathbf{E} \bigl[ \overline{v}_{n+1}(Z_{n+1}) \mathbf{1}_{\{S_{n+1}
\leq
r_{n+1,\beta}(Z_{n},S_n)\}} | Z_{n} \bigr] \\
&= & \mathbf{1}_{ \{ r_{n+1,\beta}(Z_{n},S_n)\geq t^{*}(Z_{n}) \} }
K\overline{v}_{n+1}(Z_{n})\\
&&{}+ \mathbf{1}_{ \{ r_{n+1,\beta}(Z_{n},S_n) < t^{*}(Z_{n}) \} }
J(\overline{v}_{n+1},g)(Z_{n},r_{n+1,\beta}(Z_{n},S_n)).
\end{eqnarray*}
However, our definition of $r_{n,\beta}$ with the special use of
parameter $\beta$ implies
\[
\{ r_{n+1,\beta}(Z_{n},S_n) \geq t^{*}(Z_{n}) \}
= \Bigl\{ \widehat{K}_{n+1}\widehat{v}_{n+1}(\widehat{Z}_{n})
> \max_{s\in G(\widehat{Z}_{n})} \widehat{J}_{n+1}(\widehat
{v}_{n+1},g)(\widehat{Z}_{n},s) \Bigr\}.
\]
Consequently, one obtains
%
%
\begin{eqnarray} \label{NF1}
\overline{v}_{n} (Z_{n}) &=& \widehat{K}_{n+1}\widehat
{v}_{n+1}(\widehat
{Z}_{n}) \vee
\max_{s\in G(\widehat{Z}_{n})} \widehat{J}_{n+1}(\widehat
{v}_{n+1},g)(\widehat{Z}_{n},s)
\nonumber\\
&&{} + \mathbf{1}_{ \{ r_{n+1,\beta}(Z_{n},S_n)\geq t^{*}(Z_{n}) \} }
[ K\overline{v}_{n+1}(Z_{n}) - \widehat{K}_{n+1}\widehat
{v}_{n+1}(\widehat{Z}_{n}) ] \nonumber\\[-8pt]\\[-8pt]
&&{} + \mathbf{1}_{ \{ r_{n+1,\beta}(Z_{n},S_n) < t^{*}(Z_{n}) \} }
\Bigl[ J(\overline{v}_{n+1},g)(Z_{n},r_{n+1,\beta}(Z_{n},S_n))\nonumber\\
&&\hspace*{109pt}{} - \max
_{s\in G(\widehat{Z}_{n})} \widehat{J}_{n+1}(\widehat
{v}_{n+1},g)(\widehat{Z}_{n},s) \Bigr].\nonumber
\end{eqnarray}
Let us study the term with operator $K$. First, we insert $V_n$ to be
able to use our work from the previous section (we cannot directly
apply it to $\overline{v}_n$ because it may not be
Lipschitz-continuous). Clearly, one has
%
%
\begin{eqnarray}\label{NF2}
&&| K\overline{v}_{n+1}(Z_{n}) - \widehat{K}_{n+1}\widehat
{v}_{n+1}(\widehat{Z}_{n}) |
\nonumber\\[-8pt]\\[-8pt]
&&\qquad \leq\mathbf{E} [|\overline{V}_{n+1}-V_{n+1}| |{Z}_{n}]
+ | Kv_{n+1}(Z_{n})-\widehat{K}_{n+1}\widehat
{v}_{n+1}(\widehat
{Z}_{n}) |.\nonumber
\end{eqnarray}
Similar calculations to those of Lemmas \ref{Tech4}, \ref{lemKquantif}
and \ref{lemerreur4}, and equation (\ref{Ldchap2solved}) yield
%
%
\begin{eqnarray}\label{NF3}
&&| Kv_{n+1}(Z_{n})-\widehat{K}_{n+1}\widehat
{v}_{n+1}(\widehat{Z}_{n}) |\nonumber\\
&&\qquad \leq (C_gE_4+ [ v_{n+1} ]_{1} E_{2}+ [ v_{n+1} ]_{*} [ Q ]
) ( | {Z}_{n}-\widehat{Z}_{n} | +\mathbf{E}
[|{Z}_{n}-\widehat{Z}_{n} | |\widehat{Z}_{n} ] )\hspace*{-27pt}\\
&&\qquad\quad{} +2[ v_{n+1} ]\mathbf{E} [|{Z}_{n+1}-\widehat{Z}_{n+1}|
|\widehat{Z}_{n} ]+ \mathbf{E} [|{V}_{n+1}-\widehat{V}_{n+1} |
|\widehat{Z}_{n} ].\nonumber
\end{eqnarray}
Now we turn to operator $J$. Set $R_n=r_{n+1,\beta}(Z_{n},S_n)$. We
first study the case when $ R_n = s^{*}_{n+1}(Z_{n},S_n) <
t^{*}(Z_{n})$. By definition, one has
\[
\widehat{J}_{n+1}(\widehat{v}_{n+1},g)(\widehat{Z}_{n},R_n) = \max
_{s\in G(\widehat{Z}_{n})} \widehat{J}_{n+1}(\widehat
{v}_{n+1},g)(\widehat{Z}_{n},s).
\]
As above, we insert $V_n$ and obtain
%
%
\begin{eqnarray}\label{NF4}
\hspace*{22pt}&&\Bigl| \Bigl[
J (\overline{v}_{n+1} ,g) (Z_{n},R_n) - \max_{s\in G(\widehat{Z}_{n})}
\widehat{J}_{n+1}(\widehat{v}_{n+1},g)(\widehat{Z}_{n},s) \Bigr]
\mathbf
{1}_{\{ R_n = s^{*}_{n+1}(Z_{n},S_n)\} } \Bigr| \nonumber\\
\hspace*{22pt}&&\qquad\leq \mathbf{E} [|\overline{V}_{n+1}-V_{n+1}|
|{Z}_{n}
]\mathbf{1}_{\{ R_n = s^{*}_{n+1}(Z_{n},S_n)\} }\\
\hspace*{22pt}&&\qquad\quad{} + |
J ({v}_{n+1},g) (Z_{n},R_n) - \widehat{J}_{n+1}(\widehat
{v}_{n+1},g)(\widehat{Z}_{n},R_n)) | \mathbf{1}_{\{ R_n =
s^{*}_{n+1}(Z_{n},S_n)\} } .\nonumber\hspace{-9pt}
\end{eqnarray}
Again, similar arguments as those used for Lemmas \ref{Tech3}, \ref
{lem_indic4} and \ref{lemerreur4}, and equations (\ref{Ld1b}),
(\ref
{Ld1bsolved}) and (\ref{Ld1asolved}) yield, on $\{ R_n =
s^{*}_{n+1}(Z_{n},S_n)\}$
%
%
\begin{eqnarray}\label{NF5}\qquad
&&
|J ({v}_{n+1},g) (Z_{n},R_n) - \widehat{J}_{n+1}(\widehat
{v}_{n+1},g)(\widehat{Z}_{n},R_n) |\nonumber\\
&&\qquad\leq \bigl([ v_{n+1} ]_{1} E_{2}+ [ g ]_{1}+C_gC_{t^*}[ \lambda
]_1(2+C_{t^*}C_{\lambda})\bigr)\nonumber\\
&&\qquad\quad{}\times ( | {Z}_{n}-\widehat{Z}_{n} |
+\mathbf{E} [|{Z}_{n}-\widehat{Z}_{n} | |\widehat{Z}_{n} ]
)\nonumber\\[-8pt]\\[-8pt]
&&\qquad\quad{} +2[ v_{n+1} ]\mathbf{E} [|{Z}_{n+1}-\widehat{Z}_{n+1}|
|\widehat{Z}_{n} ]+ \mathbf{E} [|{V}_{n+1}-\widehat{V}_{n+1} |
|\widehat{Z}_{n} ]\nonumber\\
&&\qquad\quad{} +[ g ]_1\mathbf{E} [|{Z}_{n}-\widehat{Z}_{n}| |\widehat{Z}_{n}]\nonumber\\
&&\qquad\quad{} + 2C_g\mathbf{E} \bigl[\bigl|\mathbf{1}_{\{{S}_{n+1}<R_n\}}-\mathbf{1}_{\{
\widehat{S}_{n+1}<R_n\}}\bigr| | \widehat{Z}_{n} \bigr].\nonumber
\end{eqnarray}
Note that all the constants with a factor $[ t^* ]$ have vanished
because we know here that both $R_n<t^{*}(Z_n)$ and $R_n<t^{*}(\widehat
{Z}_n)$ hold on $\{ R_n = s^{*}_{n+1}(Z_{n},S_n)\}$.

Finally, on $\{s^{*}(Z_{n}) \geq t^{*}(Z_{n})=R_n+ \beta\}$,
by construction of the grid $G(\widehat{Z}_{n})$ (see Remark \ref
{rqgrilletstar}), one has for all $0<\eta<\min_{z\in\Gamma
^Z_{n}}\{
\Delta(z)\}$,
\[
R_n = t^{*}(Z_{n}) - \beta< s^{*}(Z_{n}) < t^{*}(\widehat
{Z}_{n})-\eta.
\]
Consequently, using the crude bound
\[
|J (\overline{v}_{n+1} ,g) (Z_{n},R_n) |
+ \Bigl| \max_{s\in G(\widehat{Z}_{n})} \widehat{J}_{n+1}(\widehat
{v}_{n+1},g)(\widehat{Z}_{n},s) \Bigr| \leq2 C_{g},
\]
one obtains
%
%
\begin{eqnarray}\label{NF6}
&&
\Bigl|J (\overline{v}_{n+1}, g) (Z_{n},r_{n+1,\beta}(Z_{n}, S_n)) - \max
_{s\in G(\widehat{Z}_{n})} \widehat{J}_{n+1}(\widehat
{v}_{n+1},g)(\widehat{Z}_{n},s) \Bigr|\nonumber\\
&&\quad{}\times \mathbf{1}_{\{ r_{n+1,\beta
}(Z_{n}, S_n) = t^{*}(Z_{n})-\beta\} } \\
&&\qquad \leq2 C_{g} \bigl| \mathbf{1}_{\{ t^{*}(Z_{n}) - \beta<
t^{*}(\widehat{Z}_{n})- \eta\} } \bigr|.\nonumber
\end{eqnarray}
Now the combination of equations (\ref{NF1})--(\ref{NF6}) and
Lemmas \ref
{lem_indic} and \ref{lem_indic4} yields, for all $\beta<\eta
<\min
_{z\in\Gamma^Z_{n}}\{\Delta(z)\}$
\begin{eqnarray*}
\|\overline{V}_n-\widehat{V}_n \|_p&\leq& \|\overline
{V}_{n+1}-{V}_{n+1} \|_p+ \|{V}_{n+1}-\widehat{V}_{n+1} \|
_p+2[ v_{n+1} ] \|Z_{n+1}-\widehat{Z}_{n+1} \|_p\\
&&{}+ \|Z_{n}-\widehat{Z}_{n} \|_p \bigl(2[ v_{n+1} ]_1E_2+2C_gC_{t^*}[
\lambda ]_1(2+C_{t^*}C_{\lambda})\\
&&\hspace*{69.4pt}{}+ (4C_gC_{\lambda}[ t^* ]+2[ v_{n+1} ]_*[ Q ] )\vee
(3[ g ]_1) \bigr)\\
&&{} +2C_g \biggl(2C_{\lambda}\eta+\frac{1}{\eta} \|S_{n+1}-\widehat{S}_{n+1}
\|_p+\frac{[ t^* ]}{\eta-\beta} \|Z_{n}-\widehat{Z}_{n} \|
_p \biggr).
\end{eqnarray*}
Now suppose there exists $0<a<1$ such that $\eta=a^{-1}\beta$. Then the
optimal choice for $\eta$ satisfies
\[
2C_{\lambda}\eta=\frac{1}{\eta} \biggl(\frac{[ t^* ]}{1-a}\|\widehat
{Z}_{n}-{Z}_{n}\|_p+\|{S}_{n+1}-\widehat{S}_{n+1}\|_p \biggr),
\]
providing it also satisfies the condition $0<\eta<\min_{z\in
\Gamma
^{z}_{n}}\{\Delta(z)\}$, hence the result.
\end{pf}

Theorem \ref{thtpsarret} gives a recursive error estimation. Here is
the initializing step.
\begin{theorem}\label{thtpsarret2}
Suppose the discretization parameters are chosen such that there exists
$0<a<1$ satisfying
\[
\frac{\beta}{a}=(2C_{\lambda})^{-1/2} \biggl(\frac{[ t^* ]}{1-a}\|
\widehat{Z}_{N-1}-{Z}_{N-1}\|_p+\|{S}_{N}-\widehat{S}_{N}\|_p \biggr)^{1/2}
<\min
_{z\in
\Gamma^{z}_{N-1}}\{\Delta(z)\}.
\]
Then one has
\begin{eqnarray*}
&&\| \overline{V}_{N-1} - V_{N-1} \|_{p} \\
&&\qquad \leq
\|\widehat{V}_{N-1} - V_{N-1} \|_{p}+3[ g ] \| {Z}_{N}-\widehat{Z}_{N}
\|_p
+a_{N-1} \| {Z}_{N-1}-\widehat{Z}_{N-1} \|_p\\
&&\qquad\quad{} +4C_g(2C_{\lambda})^{1/2} \biggl(\frac{[ t^* ]}{1-a}\|\widehat
{Z}_{N-1}-{Z}_{N-1}\|_p+\|{S}_{N}-\widehat{S}_{N}\|_p \biggr)^{1/2}
\end{eqnarray*}
with $a_{N-1}= (2[ g ]_1E_2+2C_gC_{t^*}[ \lambda
]_1(2+C_{t^*}C_{\lambda})+ (4C_gC_{\lambda}[ t^* ]+2[ g ]_*[ Q ]
)\vee(3[ g ]_1) )$.
\end{theorem}
\begin{pf}
As before, the strong Markov property of the process $\{X(t)\}$ yields
\begin{eqnarray*}
\overline{v}_{N-1}(Z_{N-1}) & = & \mathbf{E} \bigl[ g\bigl(X_{r_{N,\beta
}(Z_{N-1}, S_{N-1})}\bigr) \mathbf{1}_{\{S_{N}> r_{N,\beta}(Z_{N-1},
S_{N-1})\}} | Z_{N-1} \bigr]\\
&&{} +\mathbf{E} \bigl[ g(Z_{N}) \mathbf{1}_{\{S_{N} \leq r_{N,\beta
}(Z_{N-1}, S_{N-1})\}} | Z_{N-1} \bigr] \nonumber\\
& = & \mathbf{1}_{ \{ r_{N,\beta}(Z_{N-1}, S_{N-1})\geq t^{*}(Z_{N-1})
\} } Kg(Z_{N-1}) \nonumber\\
&&{} + \mathbf{1}_{ \{ r_{N,\beta}(Z_{N-1}, S_{N-1}) < t^{*}(Z_{N-1}) \} }
J(g,g)(Z_{N-1},r_{N,\beta}(Z_{N-1}, S_{N-1})).
\end{eqnarray*}
The rest of the proof is similar to that of the previous theorem.
\end{pf}

As in Section \ref{section_erreur}, it is now clear that an adequate
choice of discretization parameters yields arbitrarily small errors if
one uses the stopping-time $\tau_N$.

\section{Example}
\label{section_appli}
Now we apply the procedures described in Sections \ref{section_Vchap}
and \ref{sectiontpsarret} on a simple PDMP and present numerical results.

\begin{figure}

\includegraphics{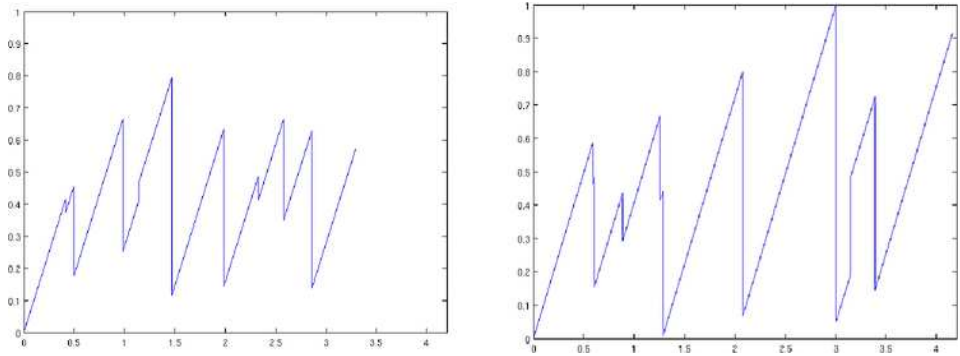}

\caption{Two trajectories of the PDMP.}\label{figure_traj}
\end{figure}

Set $E=[0, 1[$ and $\partial E=\{1\}$. The flow is defined on $[0,1]$
by $\phi(x,t)=x+vt$ for some positive $v$, the jump rate is defined on
$[0,1]$ by $\lambda(x)=\beta x^{\alpha}$, with $\beta>0$ and $\alpha
\geq1$, and for all $x\in[0,1]$, one sets $Q(x,\cdot)$ to be the
uniform law on $[0, 1/2]$. Thus the process moves with constant speed
$v$ toward $1$, but the closer it gets to the boundary $1$, the higher
the probability to jump backward on $[0, 1/2]$. Figure \ref
{figure_traj} shows two trajectories of this process for $x_0=0$,
$v=\alpha=1$ and $\beta=3$ and up to the $10$th jump.

The reward function $g$ is defined on $[0,1]$ by $g(x)=x$. Our
assumptions are clearly satisfied, and we are even in the special case
when the flow is Lipschitz-continuous (see Remark \ref{remarquelipvn}).
All the constants involved in Theorems \ref{th_Vn} and \ref{thtpsarret}
can be computed explicitly.

The real value function $V_0=v_0(x_0)$ is unknown, but, as our stopping
rule $\tau_N$ is a stopping time dominated by $T_N$, one clearly has
%
%
\begin{eqnarray}\label{empiric}
\overline{V}_0&=&\mathbf{E}_{x_{0}} [ g (X({\tau_{N}})
) ]
\leq V_0=\sup_{\tau\in\mathcal{M}_{N}}\mathbb{E}_{x_{0}}
[g(X(\tau
)) ]\nonumber\\[-8pt]\\[-8pt]
&\leq&\mathbb{E}_{x_{0}} \Bigl[\sup_{0\leq t\leq T_N}g (X(t)
) \Bigr].\nonumber
\end{eqnarray}
%
%
%
The first and last terms can be evaluated by Monte Carlo simulations,
which provide another indicator of the sharpness of our numerical
procedure. For $10^6$ Monte Carlo simulations, one obtains $\mathbb
{E}_{x_{0}} [\sup_{0\leq t\leq T_N}g (X(t) ) ]=0.9878$.
Simulation results (for $d=2$, $x_0=0$, $v=\alpha=1$, $\beta=3$, up to
the $10$th jump and for $10^5$ Monte Carlo simulations) are given in
Table \ref{table_results}. Note that, as expected, the theoretical
%
%
%
\begin{table}
\caption{Simulation results}\label{table_results}
\begin{tabular*}{\tablewidth}{@{\extracolsep{\fill}}lccccccc@{}}
\hline
$\bolds{Pt}$ & $\bolds{QE}$ & $\bolds\Delta$
& $\bolds{\widehat{V}_0}$ & $\bolds{\overline{V}_0}$
& $\bolds{B_{1}}$ & $\bolds{B_{2}}$ & $\bolds{B_{3}}$
\\
\hline
\phantom{0}10&0.0943&0.151&0.7760&0.8173&0.1705&74.64&897.0\\
\phantom{0}50&0.0418&0.100&0.8298&0.8785&0.1093&43.36&511.5\\
100&0.0289&0.083&0.8242&0.8850&0.1028&34.15&400.3\\
500&0.0133&0.056&0.8432&0.8899&0.0989&21.03&243.1\\
900&0.0102&0.049&0.8514&0.8968&0.0910&17.98&206.9\\
\hline
\end{tabular*}
\tablewidth=230pt
\begin{tabular*}{\tablewidth}{@{\extracolsep{\fill}}ll@{}}
\hline
$Pt$&Number of points in each quantization grid\\
$QE$&Quantization error: $QE={\max_{0\leq k\leq N}}\|\Theta_k-\widehat
{\Theta }_k\|_2$\\
$\Delta$& For all $z$, $\Delta(z)=\Delta$\\
$B_{1}$&Empirical bound $\mathbb{E}_{x_{0}} [\sup_{0\leq t\leq
T_N}g (X(t) ) ]-\overline{V}_0$\\
$B_{2}$&Theoretical bound given by Theorem \ref{th_Vn}\\
$B_{3}$&Theoretical bound given by Theorems \ref{thtpsarret} and
\ref{thtpsarret2}\\
\hline
\end{tabular*}
\end{table}
errors decrease as the quantization error decreases.
From (\ref{empiric}), it follows that
\[
V_0-\overline{V}_0 \leq\mathbb{E}_{x_{0}} \Bigl[\sup_{0\leq t\leq
T_N}g (X(t) ) \Bigr]-\overline{V}_0.
\]
This provides an empirical upper bound for the error.

%
\begin{appendix}\label{app}
\section{Auxiliary results}\label{app_A}
\subsection{Lipschitz properties of $J$ and $K$}
In this section, we derive useful Lipschitz-type properties of
operators $J$ and $K$. The first result is straightforward.
\begin{lemma}
\label{Tech1}
Let $h\in\mathbf{L}^{c}$. Then for all $(x,y)\in E^{2}$ and $(t,u)\in
\mathbb{R}_+^2$, one has
\begin{eqnarray*}
&&
\bigl| h \bigl( \phi\bigl(x,t\wedge t^{*}(x)\bigr) \bigr) {e}^{-\Lambda
(x,t\wedge t^{*}(x))} - h \bigl(\phi\bigl(y,u\wedge t^{*}(y)\bigr) \bigr)
{e}^{-\Lambda(y,u\wedge t^{*}(y))} \bigr| \\
&&\qquad \leq D_{1}(h) |x-y| +D_{2}(h) |t-u|,
\end{eqnarray*}
where:
\begin{itemize}
\item if $t<t^{*}(x)$ and $u<t^{*}(y)$,
\[
D_1(h)=[ h ]_{1}+C_{h}C_{t^{*}}[ \lambda ]_{1},\qquad D_2(h)=[ h ]_2+C_{h}
C_{\lambda},
\]
\item if $t=t^{*}(x)$ and $u=t^{*}(y)$,
\[
D_1(h)=[ h ]_*+C_{h}C_{t^{*}}[ \lambda ]_1+C_{h}C_{\lambda}[ t^{*} ],
\qquad
D_{2}(h)=0,
\]
\item otherwise,
\[
D_1(h)=[ h ]_{1}+C_{h}C_{t^{*}}[ \lambda ]_1+[ h ]_{2} [ t^{*}
]+C_{h}C_{\lambda} [ t^{*} ],\qquad D_2(h)=[ h ]_{2}+C_{h}
C_{\lambda}.
\]
\end{itemize}
\end{lemma}
%
%
\begin{lemma}
\label{Tech2}
Let $w\in\mathbf{B}(E)$. Then for all $x\in E$, $(t,u)\in\mathbb{R}
_{+}^{2}$, one has
\[
| J(w,g)(x,t) - J(w,g)(x,u)  | \leq( C_{w}
C_{\lambda
} + [ g ]_{2} +C_{g} C_{\lambda} ) |t-u|.
\]
\end{lemma}
\begin{pf}
By definition of $J$, we obtain
\begin{eqnarray*}
&&
| J(w,g)(x, t) - J(w,g)(x,u)  | \\
&&\qquad\leq\biggl| \int
_{t\wedge t^{*}(x)}^{u\wedge t^{*}(x)} \lambda Qw(\phi(x,s))
e^{-\Lambda
(x,s)} \,ds \biggr| \\
&&\qquad\quad{} + \bigl| g\bigl(\phi\bigl(x,t\wedge t^{*}(x)\bigr)\bigr) e^{-\Lambda(x,t\wedge t^{*}(x))}
- g\bigl(\phi\bigl(x,u\wedge t^{*}(x)\bigr)\bigr) e^{-\Lambda(x,u\wedge t^{*}(x))}\bigr|.
\end{eqnarray*}
Applying Lemma \ref{Tech1} to $h=g$, the result follows.
\end{pf}
\begin{lemma}
\label{Tech3}
Let $w\in\mathbf{L}^{c}$. Then for all $(x,y) \in E^{2}$, $t\in
\mathbb{R}_{+}$,
\[
| J(w,g)(x,t) - J(w,g)(y,t)  | \leq( C_{w} E_{1}
+ [ w ]_{1} E_{2} + E_{3} ) |x-y|,
\]
where
\begin{eqnarray*}
E_{1} & = & C_{\lambda} [ t^{*} ] + C_{t^{*}} [ \lambda ]_{1}
(
1+ C_{t^{*}} C_{\lambda} ), \\
E_{2} & = & C_{t^{*}} C_{\lambda} [ Q ], \\
E_{3} & = & [ g ]_{1} + [ g ]_{2} [ t^{*} ] + C_{g} \{ C_{t^{*}}
[ \lambda ]_{1} + C_{\lambda} [ t^{*} ] \}.
\end{eqnarray*}
\end{lemma}
\begin{pf}
Again by definition, we obtain
\begin{eqnarray*}
&&
| J(w,g) (x,t) - J(w,g)(y,t)  | \nonumber\\
&&\qquad \leq\biggl| \int_{0}^{t\wedge t^{*}(x)} \lambda Qw(\phi(x,s))
e^{-\Lambda(x,s)} \,ds
- \int_{0}^{t\wedge t^{*}(y)} \lambda Qw(\phi(y,s)) e^{-\Lambda(y,s)}
\,ds \biggr| \nonumber\\
&&\qquad\quad{} + \bigl| g\bigl(\phi\bigl(x,t\wedge t^{*}(x)\bigr)\bigr) e^{-\Lambda(x,t\wedge t^{*}(x))}
- g\bigl(\phi\bigl(y,t\wedge t^{*}(y)\bigr)\bigr) e^{-\Lambda(y,t\wedge t^{*}(y))}\bigr|.
\end{eqnarray*}
Without loss of generality it can be assumed that $t^{*}(x)\leq
t^{*}(y)$. From Lemma \ref{Tech1} for $h=g$ and using the fact that
$ |t\wedge t^{*}(x)-t\wedge t^{*}(y) | \leq| t^{*}(x)-t^{*}(y) |$, we get
\begin{eqnarray*}
&&
| J(w,g) (x,t) - J(w,g)(y,t)  | \nonumber\\
&&\qquad \leq\int_{0}^{t\wedge t^{*}(x)} \bigl| \lambda Qw(\phi(x,s))
e^{-\Lambda(x,s)} - \lambda Qw(\phi(y,s)) e^{-\Lambda(y,s)} \bigr| \,ds
\nonumber\\
&&\qquad\quad{} + ( C_{w} C_{\lambda} [ t^{*} ] + E_{3} ) |x-y|.
\end{eqnarray*}
By using a similar results as Lemma \ref{Tech1} for $h=\lambda Qw$, we
obtain the result.
\end{pf}
\begin{lemma}
\label{Tech4}
Let $w\in\mathbf{L}^{c}$. Then for all $(x,y) \in E^{2}$,
\[
| Kw(x) - Kw(y) | \leq( C_{w} E_{4} + [ w ]_{1}
E_{2}+[ w ]_*[ Q ] ) |x-y|,
\]
where $E_{4} = 2C_{\lambda} [ t^{*} ] + C_{t^{*}} [ \lambda ]_{1}
( 2 + C_{t^{*}} C_{\lambda} )$.
\end{lemma}
\begin{pf}
The proof is similar to the previous ones and is therefore omitted.
\end{pf}

\subsection{Lipschitz properties of the value functions}
\label{appvnlip}
Now we turn to the Lipschitz continuity of the sequence of value
functions $(v_n)$. Namely, we prove that under our assumptions, $v_n$
belongs to $\mathbf{L}^{c}$ for all $0\leq n\leq N$. We also compute
the Lipschitz constant of $v_n$ on $\overline{E}$ as it is much sharper
in this case than $[ v_n ]_1$ (see Remark \ref{remarque_lip}).

We start with proving sharper results on operator $J$.
\begin{lemma}
\label{Tech5}
Let $w\in\mathbf{L}^{c}$. Then for all $x \in E$ and $(s,t)\in
\mathbb{R}_{+}^{2}$,
\[
\Bigl| \sup_{u\geq t}J(w,g)(x,u) - \sup_{u\geq s}J(w,g)(x,u) \Bigr|
\leq
( C_{w} C_{\lambda} + [ g ]_{2} +C_{g} C_{\lambda} )
|t-s| .
\]
\end{lemma}
\begin{pf}
Without loss of generality it can be assumed that $t\leq s$.
Therefore, one has
%
%
\begin{eqnarray}\label{eqT1}
&&
\Bigl| \sup_{u\geq t}J(w,g)(x,u) - \sup_{u\geq s}J(w,g)(x,u)
\Bigr|\nonumber\\[-8pt]\\[-8pt]
&&\qquad =
\sup_{u\geq t}J(w,g)(x,u) - \sup_{u\geq s}J(w,g)(x,u).\nonumber
\end{eqnarray}
Note that there exists $\overline{t} \in[t\wedge t^{*}(x),t^{*}(x)]$
such that
$\sup_{u\geq t}J(w,g)(x,u) = J(w,g)(x,\overline{t})$.
Consequently, if $\overline{t}\geq s$ then one has $|{ \sup
_{u\geq t}J(w,g)(x,u)} - \sup_{u\geq s}J(w,g)(x,u) |=0$.

Now if $\overline{t} \in[t\wedge t^{*}(x),s[$, then one has
\[
\sup_{u\geq t}J(w,g)(x,u) - \sup_{u\geq s}J(w,g)(x,u) \leq
J(w,g)(x,\overline{t}) - J(w,g)(x,s).
\]
From Lemma \ref{Tech2}, we obtain the following inequality:
%
%
\begin{equation}\label{eqT2}\qquad\quad
\sup_{u\geq t}J(w,g)(x,u) - \sup_{u\geq s}J(w,g)(x,u) \leq( C_{w}
C_{\lambda} + [ g ]_{2} +C_{g} C_{\lambda} ) |\overline{t}-s|.
\end{equation}
Combining (\ref{eqT1}), (\ref{eqT2}) and the fact that
$|\overline{t}-s| \leq|t-s|$ the result follows.
\end{pf}

Similarly, we obtain the following result.
\begin{lemma}
\label{Tech5bis}
Let $w\in\mathbf{L}^{c}$. Then for all $(x,y) \in E^{2}$,
\[
\Bigl| \sup_{t\leq t^{*}(x)}J(w,g)(x,t) - \sup_{t\leq t^{*}(y)}J(w,g)(y,t)
\Bigr| \leq( C_{w} {E}_{5} + [ w ]_1 E_{2} +
{E}_{6} ) |x-y|,
\]
where ${E}_{5} = E_{1} + C_{\lambda} [ t^{*} ]$ and $E_{6} = E_{3} +
( [ g ]_2+ C_{g} C_{\lambda} ) [ t^{*} ]$.
\end{lemma}

Now we turn to $(v_n)$. Recall from \cite{gugerli86} that for all
$0\leq n\leq N$, $(v_n)$ is bounded with $C_{v_n}=C_g$.
\begin{proposition}
\label{prop1}
For all $0\leq n\leq N$, $v_{n}\in\mathbf{L}^{c}$ and
%
%
\begin{eqnarray}\quad\qquad
\label{vLip1}
[ v_{n} ]_{1} &\leq& e^{C_{\lambda}C_{t^{*}}} \bigl( 2[ v_{n+1} ]_{1}
E_{2} +C_{g}E_{1}+ C_g E_4 +C_{g}C_{t^{*}}[ \lambda ]_{1}
(1+C_{\lambda}C_{t^{*}}) \bigr) \nonumber\\[-8pt]\\[-8pt]
&&{} + e^{C_{\lambda}C_{t^{*}}}
\{ ( [ g ]_1+[ g ]_2[ t^* ] )\vee( [ v_{n+1} ]_*[ Q ] ) \},\nonumber\\
\label{vLip2}
{}[ v_{n} ]_{2} & \leq& e^{C_{\lambda}C_{t^{*}}} \{
C_{g}C_{\lambda
} (4+C_{\lambda}C_{t^{*}} )
+[ g ]_{2} \},\\
{}[ v_{n} ]_*&\leq&[ v_{n} ]_{1}+[ v_{n} ]_{2}[ t^{*} ],\nonumber\\
{}[ v_{n} ]&\leq&[ v_{n+1} ]_1 E_2+ C_g E_5 + \{E_6\vee
([ v_{n+1} ]_*[ Q ] +C_{g}C_{t^{*}}[ \lambda ]_{1} ) \}
.\nonumber
\end{eqnarray}
\end{proposition}
\begin{pf}
Clearly, $v_{N}=g$ is in $\mathbf{L}^{c}$. Assume that $v_{n+1}$ is in
$\mathbf{L}^{c}$, then
by using the semi-group property of the drift $\phi$ it can be shown
that for any $x\in E$, $t\in[0,t^{*}(x)]$, one has (see \cite
{gugerli86}, equation (8))
%
%
\begin{eqnarray}\label{vlipt1}
&&v_{n}(\phi(x,t)) = e^{\Lambda(x,t)} \Bigl\{ \Bigl( \sup_{u\geq
t}J(v_{n+1},g)(x,u) \vee Kv_{n+1}(x) \Bigr)\nonumber\\[-8pt]\\[-8pt]
&&\hspace*{191.2pt}{} -Iv_{n+1}(x,t)
\Bigr\}.\nonumber
\end{eqnarray}
Note that for $x\in E$, $t\in\mathbb{R}_{+}$, one has
%
%
\begin{eqnarray}\label{vlipt1b}
&&\sup_{u\geq t}J(v_{n+1},g)(x,u) \vee Kv_{n+1}(x) \nonumber\\
&&\qquad\leq \sup
_{u}J(v_{n+1},g)(x,u) \vee Kv_{n+1}(x)\\
&&\qquad=v_{n}(x).\nonumber
\end{eqnarray}
Set $(x,y)\in E^{2}$ and $t\in[0,t^{*}(x)\wedge t^{*}(y)]$. It is
easy to show that
%
%
\begin{eqnarray}
\bigl| e^{\Lambda(x,t)}-e^{\Lambda(y,t)} \bigr| & \leq&
e^{C_{\lambda}
C_{t^{*}}} [ \lambda ]_{1} C_{t^{*}} |x-y|, \\
\label{vlipt2}
| Iv_{n+1}(x,t)-Iv_{n+1}(y,t) | & \leq& ( C_{v_{n+1}}
E_{1} + [ v_{n+1} ]_{1} E_{2} ) |x-y|.
\end{eqnarray}
Then, (\ref{vlipt1})--(\ref{vlipt2}) yield
%
%
\begin{eqnarray}\label{vlipt2b}
&&
|v_{n} (\phi(x,t))-v_{n}(\phi(y,t)) | \nonumber\\
&&\qquad
\leq\{ |
v_{n}(x) | + | Iv_{n+1}(x,t) | \} e^{C_{\lambda}
C_{t^{*}}} [ \lambda ]_{1} C_{t^{*}} |x-y| \nonumber\\
&&\qquad\quad{} +e^{\Lambda(y,t)} \Bigl\{ \sup_{u\geq t} | J(v_{n+1},g)(x,u) -
J(v_{n+1},g)(y,u) | \\
&&\qquad\quad\hspace*{108.3pt}{}\vee
| Kv_{n+1}(x) - Kv_{n+1}(y) | \Bigr\} \nonumber\\
&&\qquad\quad{} +e^{\Lambda(y,t)} ( C_{v_{n+1}} E_{1} + [ v_{n+1} ]_{1} E_{2}
) |x-y|.\nonumber
\end{eqnarray}
For $x\in E$, $t\in[0,t^{*}(x)]$ and $n\in\mathbb{N}$,
note that
%
%
\begin{eqnarray}\label{vlipt2bb}
e^{\Lambda(x,t)}&\leq& e^{C_{\lambda} C_{t^{*}}},\nonumber\\[-8pt]\\[-8pt]
|Iv_{n+1}(x,t) |&\leq& C_{\lambda} C_{v_{n+1}} C_{t^{*}}\quad
\mbox{and}\quad | v_{n+1}(x) | \leq C_{g}.\nonumber
\end{eqnarray}
Therefore, we obtain inequality (\ref{vLip1}) by using (\ref
{vlipt2b}), (\ref{vlipt2bb}) and Lemma~\ref{Tech3}, \ref{Tech5}, and
the fact that $C_g E_1+E_3= C_g E_4+[ g ]_1+[ g ]_2[ t^* ]$.

Now, set $x\in E$ and $t$, $s\in[0,t^{*}(x)]$. Similarly,
one has
%
%
\begin{eqnarray}
\label{vlipt3}
\bigl| e^{\Lambda(x,t)}-e^{\Lambda(x,s)} \bigr| & \leq&
e^{C_{\lambda}
C_{t^{*}}} C_{\lambda} |t-s|,
\\
\label{vlipt4}
| Iv_{n+1}(x,t)-Iv_{n+1}(x,s) | & \leq& C_{\lambda}
C_{v_{n+1}} |t-s|.
\end{eqnarray}
Combining (\ref{vlipt1}), (\ref{vlipt1b}), (\ref{vlipt3}) and
(\ref{vlipt4}), it yields
%
%
\begin{eqnarray}\label{vlipt5}\quad
&&| v_{n}  (\phi(x,t))-v_{n}(\phi(x,s)) | \nonumber\\
&&\qquad\leq\{ |
v_{n}(x) | + | Iv_{n+1}(x,t) | \} e^{C_{\lambda}
C_{t^{*}}} C_{\lambda} |t-s| \nonumber\\[-8pt]\\[-8pt]
&&\qquad\quad{} +e^{\Lambda(x,t)} \Bigl\{ \Bigl| \sup_{u\geq t}J(v_{n+1},g)(x,u) -
\sup
_{u\geq s}J(v_{n+1},g)(x,u) \Bigr|\nonumber\\
&&\qquad\quad\hspace*{168.2pt}{} + C_{\lambda} C_{v_{n+1}}
|t-s|\Bigr\}.\nonumber
\end{eqnarray}
Finally, inequality (\ref{vLip2}) follows from equations (\ref
{vlipt2bb}), (\ref{vlipt5}) and Lemma \ref{Tech4}.

One clearly has $[ v_{n} ]_*\leq[ v_{n} ]_{1}+[ v_{n} ]_{2}[ t^{*} ]$.
Finally, set $(x,y)\in\overline{E}{}^{2}$. By
definition, one has
\begin{eqnarray*}
&&|v_n(x)-v_n(y) |\\
&&\qquad\leq \Bigl| \sup_{u\leq t^{*}(x)}J(v_{n+1},g)(x,u) - \sup_{u\leq
t^{*}(y)}J(v_{n+1},g)(y,u) \Bigr|\\
&&\qquad\quad{}\vee| Kv_{n+1}(x) - Kv_{n+1}(y)|
\end{eqnarray*}
and we conclude using Lemmas \ref{Tech5bis} and \ref{Tech4}, and the
fact that $E_4=E_5+C_{t^*}[ \lambda ]_1$.
\end{pf}
\begin{remark}\label{remarquelipvn}
Note that $[ v_n ]$ is much sharper than $[ v_n ]_1$. If in addition
to our assumptions, the drift $\phi$ is Lipschitz-continuous in both
variables, then with obvious notation, one has $[ v_n ]_i\leq[ v_n ][
\phi ]_i$ for $i\in\{1,2,*\}$, which should yield better
constants (see, e.g., Section \ref{section_appli}).
\end{remark}

\section{Structure of the stopping times of PDMP\textup{s}}
\label{app_B}
Let $\tau$ be an $\{\mathcal{F}_{t}\}_{t\in\mathbb{R}
_{+}}$-stopping time. Let us recall the important result from Davis
\cite{davis93}.
\begin{theorem}
There exists a sequence of nonnegative random variables $
(R_{n})_{n\in\mathbb{N}^{*}}$ such that $R_{n}$ is
$\mathcal{F}_{T_{n-1}}$-measurable and $\tau\wedge T_{n+1}=
(T_{n}+R_{n+1})\wedge T_{n+1}$ on
$\{\tau\geq T_{n}\}$.
\end{theorem}
\begin{lemma}
\label{lem1}
Define $\overline{R}_{1}=R_{1}$, and $\overline{R}_{k}=R_{k} \mathbf
{1}_{\{S_{k-1}\leq\overline{R}_{k-1}\}}$. Then one has
\[
\tau= \sum_{n=1}^{\infty} \overline{R}_{n}\wedge S_{n}.
\]
\end{lemma}
\begin{pf}
Clearly, on $\{T_{k}\leq\tau< T_{k+1}\}$, one has $R_{j}\geq S_{j}$
and $R_{k+1} < S_{k+1}$ for all $j\leq k$.
Consequently, by definition $\overline{R}_{j}=R_{j}$ for all $j\leq
k+1$, whence
\begin{eqnarray*}
\sum_{n=1}^{\infty} \overline{R}_{n}\wedge S_{n} & = & \sum_{n=1}^{k}
\overline{R}_{n}\wedge S_{n} + \{ \overline{R}_{k+1}\wedge
S_{k+1} \} + \sum_{n=k+2}^{\infty} \overline{R}_{n}\wedge S_{n}
\\
& = & T_{k} + R_{k+1} + \sum_{n=k+2}^{\infty} \overline{R}_{n}\wedge S_{n}.
\end{eqnarray*}
Since $\overline{R}_{k+1}=R_{k+1}< S_{k+1}$ we have $\overline
{R}_{j}=0$ for all $j\geq k+2$.
Therefore, $\sum_{n=1}^{\infty} \overline{R}_{n}\wedge S_{n}= T_{k} +
R_{k+1}= \tau$, showing the result.
\end{pf}

There exists a sequence of measurable mappings $(r_{k}
)_{k\in
\mathbb{N}_{*}}$ defined on $E\times(\mathbb{R}_{+}\times E)^{k-1}$
with value in
$\mathbb{R}_{+}$ satisfying
\begin{eqnarray*}
R_{1} & = & r_{1}(Z_{0}), \\
R_{k} & = & r_{k}(Z_{0},\Gamma_{k-1}),
\end{eqnarray*}
where $\Gamma_{k}= (S_{1},Z_{1},\ldots,S_{k},Z_{k} )$.
\begin{definition}
\label{defgene}
Consider $p\in\mathbb{N}_{*}$.
Let $(\widehat{R}_{k})_{k\in\mathbb{N}_{*}}$ be a sequence of
mappings defined on
$E\times(\mathbb{R}_{+}\times E)^{p} \times\Omega$ with value in
$\mathbb{R}_{+}$
defined by
\[
\widehat{R}_{1}(y,\gamma,\omega) = r_{p+1}(y,\gamma)
\]
and for $k\geq2$
\[
\widehat{R}_{k}(y,\gamma,\omega) = r_{p+k}(y,\gamma,\Gamma
_{k-1}(\omega)) \mathbf{1}_{\{S_{k-1}\leq\widehat{R}_{k-1}\}}
(y,\gamma
,\omega).
\]
\end{definition}
\begin{proposition}
\label{Aprop1}
Assume that $T_{p} \leq\tau\leq T_{N}$.
Then, one has
\[
\tau= T_{p} + \widehat{\tau}(Z_{0},\Gamma_{p},\theta_{T_{p}}),
\]
where $\widehat{\tau} \dvtx E\times(\mathbb{R}_{+}\times E)^{p}
\times
\Omega\to\mathbb{R}_{+}$ is defined by
%
%
\begin{equation}\label{defTauhat}
\widehat{\tau}(y,\gamma,\omega) = \sum_{n=1}^{N-p} \widehat
{R}_{n}(y,\gamma,\omega) \wedge S_{n}(\omega).
\end{equation}
\end{proposition}
\begin{pf}
First, let us prove by induction that for $k\in\mathbb{N}_{*}$, one has
%
%
\begin{equation}\label{ind}
\widehat{R}_{k}(Z_{0},\Gamma_{p},\theta_{T_{p}}) = \overline{R}_{p+k}.
\end{equation}
Indeed, one has $\widehat{R}_{1}(Z_{0},\Gamma_{p},\theta_{T_{p}}) = R_{p+1}$,
and on the set $\{\tau\geq T_{p}\}$, one also has $R_{p+1}=\overline
{R}_{p+1}$.
Consequently, $\widehat{R}_{1}(Z_{0},\Gamma_{p}) =\overline{R}_{p+1}$.
Now assume that $\widehat{R}_{k}(Z_{0},\Gamma_{p}$,\break $\theta_{T_{p}}) =
\overline{R}_{p+k}$.
Then, one has
\begin{eqnarray*}
&&
\widehat{R}_{k+1}(Z_{0}(\omega),\Gamma_{p}(\omega),\theta
_{T_{p}}(\omega))
\\
&&\qquad= r_{p+k+1}(Z_{0}(\omega), \Gamma_{p}(\omega),\Gamma_{k}(\theta
_{T_{p}}(\omega)))
\mathbf{1}_{\{S_{k}\leq\widehat{R}_{k}\}} (Z_{0}(\omega),\Gamma
_{p}(\omega),\theta_{T_{p}}(\omega)).
\end{eqnarray*}
By definition, one has $\Gamma_{k}(\theta_{T_{p}}(\omega))=
(S_{p+1}(\omega),Z_{p+1}(\omega),\ldots,S_{p+k}(\omega
),Z_{p+k}(\omega
) )$ and the induction hypothesis easily yields
$\mathbf{1}_{\{S_{k}\leq\widehat{R}_{k}\}} (Z_{0}(\omega),\Gamma
_{p}(\omega),\theta_{T_{p}}(\omega))
=\mathbf{1}_{\{S_{p+k}\leq\overline{R}_{p+k}\}} (\omega)$.
Therefore,\vspace*{1pt} we get
$\widehat{R}_{k+1}(Z_{0},\Gamma_{p},\theta_{T_{p}}) = \overline{R}_{p+k+1}$,
showing (\ref{ind}).

Combining (\ref{defTauhat}) and (\ref{ind}) yields
%
%
\begin{equation}\label{prop1eq1}
\widehat{\tau}(Z_{0},\Gamma_{p},\theta_{T_{p}}) = \sum_{n=1}^{N-n}
\overline{R}_{p+n} \wedge S_{p+n} .
\end{equation}
However, we have already seen that on the set $\{T\geq T_{p}\}$, one
has $R_{k}=\overline{R}_{k} \geq S_{k}$, for $k\leq p$.
Consequently, using (\ref{prop1eq1}), we obtain
\[
T_{p} + \widehat{\tau}(Z_{0},\Gamma_{p},\theta_{T_{p}}) = \sum
_{k=1}^{p} S_{k} + \sum_{k=p+1}^{N} \overline{R}_{k} \wedge S_{k}
=\sum_{k=1}^{N} \overline{R}_{k} \wedge S_{k}.
\]
Since $\tau\leq T_{N}$, we obtain from Lemma \ref{lem1} and its proof
that $\tau= \sum_{n=1}^{N} \overline{R}_{n}\wedge S_{n}$, showing
the result.
\end{pf}
\begin{proposition}
\label{Aprop2}
Let $(U_{n})_{n\in\mathbb{N}^{*}}$ be a sequence of nonnegative
random variables such that $U_{n}$
is $\mathcal{F}_{T_{n-1}}$-measurable and
$U_{n+1}=0$ on $\{ S_{n} > U_{n}\}$, for all $n\in\mathbb{N}_{*}$.
Set
\[
U = \sum_{n=1}^{\infty} U_{n}\wedge S_{n}.
\]
Then $U$ is an $\{\mathcal{F}_{t}\}_{t\in\mathbb{R}
_{+}}$-stopping time.
\end{proposition}
\begin{pf}
Assumption \ref{A1} yields
%
%
\begin{eqnarray}\label{prop2eq1}
\{U\leq t\} &=& \bigcup_{n=0}^{\infty} [ (\{T_{n}\leq U <
T_{n+1}\}\cap\{U\leq t\}\cap\{t<T_{n+1}\} ) \nonumber\\[-8pt]\\[-8pt]
&&\hspace*{17.6pt}{} \cup(\{T_{n}\leq U < T_{n+1}\}\cap\{U\leq t\}\cap\{
T_{n+1} \leq t \} ) ].\nonumber
\end{eqnarray}
From the definition of $U_{n}$, one has $\{U\geq T_n\}=\{U_n\geq S_n\}
$; hence one has
\begin{eqnarray*}
&&\{T_{n}\leq U < T_{n+1}\}\cap\{U\leq t\}\cap\{t<T_{n+1}\}\\
&&\qquad = \{
S_{n} \leq U_{n}\}\cap\{T_{n}+U_{n+1}\leq t\}
\cap\{T_{n} \leq t \} \cap\{t<T_{n+1} \}.
\end{eqnarray*}
Theorem 2.10(ii) in \cite{elliott82} now yields
$\{S_{n} \leq U_{n}\}\cap\{T_{n}+U_{n+1}\leq t\}\cap\{T_{n} \leq
t \} \in\mathcal{F}_{t}$; thus one has
%
%
\begin{equation}\label{prop2eq2}
\{T_{n}\leq U < T_{n+1}\}\cap\{U\leq t\}\cap\{t<T_{n+1}\}\in
\mathcal{F}_{t}.
\end{equation}
On the other hand, one has
\begin{eqnarray*}
&&
\{T_{n}\leq U < T_{n+1}\}\cap\{U\leq t\}\cap\{T_{n+1} \leq t \}\\
&&\qquad =
\{S_{n} \leq U_{n}\}\cap\{U_{n+1} < S_{n+1} \} \cap\{T_{n+1} \leq
t\}.
\end{eqnarray*}
Hence Theorem 2.10(ii) in \cite{elliott82} again yields
%
%
\begin{equation}\label{prop2eq3}
\{T_{n}\leq U < T_{n+1}\}\cap\{U\leq t\}\cap\{T_{n+1} \leq t \}
\in
\mathcal{F}_{t}.
\end{equation}
Combining equations (\ref{prop2eq1}), (\ref{prop2eq2}) and (\ref
{prop2eq3}) we obtain the result.
\end{pf}
\begin{corollary}
\label{Acoro1}
For any $(y,\gamma)\in E\times( \mathbb{R}_{+}\times E )^{p}$,
$\widehat{\tau}(y,\gamma,\cdot)$ is an $\{\mathcal
{F}_{t}\}
_{t\in\mathbb{R}_{+}}$-stopping time satisfying $\widehat{\tau
}(y,\gamma
,\cdot) \leq T_{N-p}$.
\end{corollary}
\begin{pf}
It follows\vspace*{1pt} form the definition of $\widehat{R}_{k}$ that $\widehat
{R}_{k}(y,\gamma,\omega)<S_{k}(\omega)$ implies
$\widehat{R}_{k+1}(y,\gamma,\omega)=0$ and the nonnegative random
variable $\widehat{R}_{k}(y,\gamma,\cdot)$ is
$\mathcal{F}_{T_{k-1}}$-measurable. Therefore, Proposition \ref{Aprop2}
yields that
$\widehat{\tau}(y,\gamma,\cdot)$ is an $\{\mathcal
{F}_{t}\}
_{t\in\mathbb{R}_{+}}$-stopping time.
Finally, by definition of $\widehat{\tau}$ [see (\ref
{defTauhat})], one has
$\widehat{\tau}(y,\gamma,\cdot)\leq\sum_{n=1}^{N-p}
S_{n}=T_{N-p}$ showing the result.
\end{pf}
\end{appendix}


%

%
\printaddresses


\begin{thebibliography}{17}

\bibitem{bally03}
%
\begin{barticle}[mr]
\bauthor{\bsnm{Bally},~\bfnm{Vlad}\binits{V.}} \AND
\bauthor{\bsnm{Pag{\`e}s},~\bfnm{Gilles}\binits{G.}}
(\byear{2003}).
\btitle{A quantization algorithm for solving multi-dimensional discrete-time
optimal stopping problems}.
\bjournal{Bernoulli}
\bvolume{9}
\bpages{1003--1049}.
\bid{doi={10.3150/bj/1072215199}, mr={2046816}}
\end{barticle}
%
\endbibitem

\bibitem{bally05}
%
\begin{barticle}[mr]
\bauthor{\bsnm{Bally},~\bfnm{Vlad}\binits{V.}},
\bauthor{\bsnm{Pag{\`e}s},~\bfnm{Gilles}\binits{G.}} \AND
\bauthor{\bsnm{Printems},~\bfnm{Jacques}\binits{J.}}
(\byear{2005}).
\btitle{A quantization tree method for pricing and hedging multidimensional
{A}merican options}.
\bjournal{Math. Finance}
\bvolume{15}
\bpages{119--168}.
\bid{doi={10.1111/j.0960-1627.2005.00213.x}, mr={2116799}}
\end{barticle}
%
\endbibitem

\bibitem{costa88}
%
\begin{barticle}[mr]
\bauthor{\bsnm{Costa},~\bfnm{O.~L.~V.}\binits{O.~L.~V.}} \AND
\bauthor{\bsnm{Davis},~\bfnm{M.~H.~A.}\binits{M.~H.~A.}}
(\byear{1988}).
\btitle{Approximations for optimal stopping of a piecewise-deterministic
process}.
\bjournal{Math. Control Signals Systems}
\bvolume{1}
\bpages{123--146}.
\bid{doi={10.1007/BF02551405}, mr={936330}}
\end{barticle}
%
\endbibitem

\bibitem{costa08}
%
\begin{barticle}[mr]
\bauthor{\bsnm{Costa},~\bfnm{O.~L.~V.}\binits{O.~L.~V.}} \AND
\bauthor{\bsnm{Dufour},~\bfnm{F.}\binits{F.}}
(\byear{2008}).
\btitle{Stability and ergodicity of piecewise deterministic {M}arkov
processes}.
\bjournal{SIAM J. Control Optim.}
\bvolume{47}
\bpages{1053--1077}.
\bid{doi={10.1137/060670109}, mr={2385873}}
\end{barticle}
%
\endbibitem

\bibitem{costa00}
%
\begin{barticle}[mr]
\bauthor{\bsnm{Costa},~\bfnm{O.~L.~V.}\binits{O.~L.~V.}},
\bauthor{\bsnm{Raymundo},~\bfnm{C.~A.~B.}\binits{C.~A.~B.}} \AND
\bauthor{\bsnm{Dufour},~\bfnm{F.}\binits{F.}}
(\byear{2000}).
\btitle{Optimal stopping with continuous control of piecewise deterministic
{M}arkov processes}.
\bjournal{Stochastics Stochastics Rep.}
\bvolume{70}
\bpages{41--73}.
\bid{mr={1785064}}
\end{barticle}
%
\endbibitem

\bibitem{davis93}
%
\begin{bbook}[mr]
\bauthor{\bsnm{Davis},~\bfnm{M.~H.~A.}\binits{M.~H.~A.}}
(\byear{1993}).
\btitle{Markov Models and Optimization}.
\bseries{Monographs on Statistics and Applied Probability}
\bvolume{49}.
\bpublisher{Chapman and Hall}, \baddress{London}.
\bid{mr={1283589}}
\end{bbook}
%
\endbibitem

\bibitem{dufour99}
%
\begin{barticle}[mr]
\bauthor{\bsnm{Dufour},~\bfnm{Fran{\c{c}}ois}\binits{F.}} \AND
\bauthor{\bsnm{Costa},~\bfnm{Oswaldo L.~V.}\binits{O.~L.~V.}}
(\byear{1999}).
\btitle{Stability of piecewise-deterministic {M}arkov processes}.
\bjournal{SIAM J. Control Optim.}
\bvolume{37}
\bpages{1483--1502}.
\bid{doi={10.1137/S0363012997330890}, mr={1710229}}
\end{barticle}
%
\endbibitem

\bibitem{elliott82}
%
\begin{bbook}[mr]
\bauthor{\bsnm{Elliott},~\bfnm{Robert~J.}\binits{R.~J.}}
(\byear{1982}).
\btitle{Stochastic Calculus and Applications}.
\bseries{Applications of Mathematics (New York)}
\bvolume{18}.
\bpublisher{Springer}, \baddress{New York}.
\bid{mr={678919}}
\end{bbook}
%
\endbibitem

\bibitem{gatarek91}
%
\begin{barticle}[mr]
\bauthor{\bsnm{G\c{a}tarek},~\bfnm{Dariusz}\binits{D.}}
(\byear{1991}).
\btitle{On first-order quasi-variational inequalities with integral terms}.
\bjournal{Appl. Math. Optim.}
\bvolume{24}
\bpages{85--98}.
\bid{doi={10.1007/BF01447736}, mr={1106927}}
\end{barticle}
%
\endbibitem

\bibitem{gray98}
%
\begin{barticle}[mr]
\bauthor{\bsnm{Gray},~\bfnm{Robert~M.}\binits{R.~M.}} \AND
\bauthor{\bsnm{Neuhoff},~\bfnm{David~L.}\binits{D.~L.}}
(\byear{1998}).
\btitle{Quantization}.
\bjournal{IEEE Trans. Inform. Theory}
\bvolume{44}
\bpages{2325--2383}.
\bid{doi={10.1109/18.720541}, mr={1658787}}
\end{barticle}
%
\endbibitem

\bibitem{gugerli86}
%
\begin{barticle}[mr]
\bauthor{\bsnm{Gugerli},~\bfnm{U.~S.}\binits{U.~S.}}
(\byear{1986}).
\btitle{Optimal stopping of a piecewise-deterministic {M}arkov process}.
\bjournal{Stochastics}
\bvolume{19}
\bpages{221--236}.
\bid{mr={872462}}
\end{barticle}
%
\endbibitem

\bibitem{kushner77}
%
\begin{bbook}[mr]
\bauthor{\bsnm{Kushner},~\bfnm{Harold~J.}\binits{H.~J.}}
(\byear{1977}).
\btitle{Probability Methods for Approximations in Stochastic Control
and for
Elliptic Equations}.
\bseries{Mathematics in Science and Engineering}
\bvolume{129}.
\bpublisher{Academic Press}, \baddress{New York}.
\bid{mr={0469468}}
\end{bbook}
%
\endbibitem

\bibitem{lenhart85}
%
\begin{barticle}[mr]
\bauthor{\bsnm{Lenhart},~\bfnm{Suzanne}\binits{S.}} \AND
\bauthor{\bsnm{Liao},~\bfnm{Yu~Chung}\binits{Y.~C.}}
(\byear{1985}).
\btitle{Integro-differential equations associated with optimal stopping
time of
a piecewise-deterministic process}.
\bjournal{Stochastics}
\bvolume{15}
\bpages{183--207}.
\bid{mr={869199}}
\end{barticle}
%
\endbibitem

\bibitem{pages98}
%
\begin{barticle}[mr]
\bauthor{\bsnm{Pag{\`e}s},~\bfnm{Gilles}\binits{G.}}
(\byear{1998}).
\btitle{A space quantization method for numerical integration}.
\bjournal{J. Comput. Appl. Math.}
\bvolume{89}
\bpages{1--38}.
\bid{doi={10.1016/S0377-0427(97)00190-8}, mr={1625987}}
\end{barticle}
%
\endbibitem

\bibitem{pages05}
%
\begin{barticle}[mr]
\bauthor{\bsnm{Pag{\`e}s},~\bfnm{Gilles}\binits{G.}} \AND
\bauthor{\bsnm{Pham},~\bfnm{Huy{\^e}n}\binits{H.}}
(\byear{2005}).
\btitle{Optimal quantization methods for nonlinear filtering with discrete-time
observations}.
\bjournal{Bernoulli}
\bvolume{11}
\bpages{893--932}.
\bid{doi={10.3150/bj/1130077599}, mr={2172846}}
\end{barticle}
%
\endbibitem

\bibitem{pages04b}
%
\begin{barticle}[mr]
\bauthor{\bsnm{Pag{\`e}s},~\bfnm{Gilles}\binits{G.}},
\bauthor{\bsnm{Pham},~\bfnm{Huy{\^e}n}\binits{H.}} \AND
\bauthor{\bsnm{Printems},~\bfnm{Jacques}\binits{J.}}
(\byear{2004}).
\btitle{An optimal {M}arkovian quantization algorithm for multi-dimensional
stochastic control problems}.
\bjournal{Stoch. Dyn.}
\bvolume{4}
\bpages{501--545}.
\bid{doi={10.1142/S0219493704001231}, mr={2102752}}
\end{barticle}
%
\endbibitem

\bibitem{pages04}
%
\begin{bincollection}[mr]
\bauthor{\bsnm{Pag{\`e}s},~\bfnm{Gilles}\binits{G.}},
\bauthor{\bsnm{Pham},~\bfnm{Huy{\^e}n}\binits{H.}} \AND
\bauthor{\bsnm{Printems},~\bfnm{Jacques}\binits{J.}}
(\byear{2004}).
\btitle{Optimal quantization methods and applications to numerical
problems in
finance}.
In \bbooktitle{Handbook of Computational and Numerical Methods in Finance}
\bpages{253--297}.
\bpublisher{Birkh\"auser}, \baddress{Boston, MA}.
\bid{mr={2083055}}
\end{bincollection}
%
\endbibitem

\end{thebibliography}
\end{document}